Dmytro Taranovsky
Last Update: December 31, 2018

# Ordinal Notation

**Abstract:** We introduce a framework for ordinal notation systems, present a family of strong yet simple systems, and give many examples of ordinals in these systems. While much of the material is conjectural, we include systems with conjectured strength beyond second order arithmetic (and plausibly beyond ZFC), and prove well-foundedness for some weakened versions.

**Contents:**

**Notes:**
• Claims not marked as theorems or propositions should be treated as conjectures (see Introduction).
• [OrdinalArithmetic.py](#) (absolute link) is a python module/program that implements comparison and ordinal arithmetic for the main ordinal notation system. See [2.4 One Variable C](#) subsection of the paper for a notational difference with the module. The reader is encouraged to try out the module.

# 1 Introduction

## 1.1 Introduction

Ordinals play a core role in set theory, and a good ordinal representation system for a theory gives us a qualitatively enhanced understanding of that theory. However, prior to this paper, we lacked reasonable ordinal notation systems for "strong" theories such as second order arithmetic. One can always extend a system to bigger ordinals, but the key obstacle is complexity. For example, by just intuitively counting past infinity without a plan, one will likely bog down even before reaching $\varepsilon_0$. However, a repeated historical pattern is that new ideas greatly simplify the upper reaches of the previous systems and reach new ordinals. It is just such a set of ideas that form the core of this paper. Specifically, we present a framework that it makes it much easier to define new ordinal notation systems, and define and analyze simple yet strong systems following that framework.

**Prerequisites:** The paper assumes basic understanding of ordinals, and parts of the paper assume more. A good introduction is in (Rathjen 2006).
**Note:** This paper was developed over time, and different main sections can mostly be read independently of each other. For example, the reader is free to read the definition of the main system right away.

**What this Paper Accomplishes:** We define and analyze several very strong ordinal notation systems presented as recursive relations between terms. The canonical assignment of ordinals to terms (see [1.3 Goals of Ordinal Analysis](#) below for what that means) is fully given for some weaker systems, and many examples of conjectured canonical assignment are given beyond that. Well-foundedness is proved for some of the systems. Strength lower bounds are not proved.
    Formally, unless stated (or the reader is satisfied) otherwise, claims in this paper about the notation systems should be treated as conjectures. While some

conjectures are speculative (with significant doubts always noted), in many cases our understanding is so detailed, coherent, and precise as to be very confident in its accuracy. This is not a substitute for proof, but is the best we have before the required bounds are proved.

## 1.2  Prior Work

Prior to the [Degrees of Reflection](#) system in this paper, the strongest reasonably simple ordinal notation system was the one for KP + $\Pi_3$ reflection (Rathjen 1994). (That system is referenced in the next paragraph, which the reader without requisite background is free to skip.) Rathjen also proposed stronger systems, but they are so complex that the problem of finding a reasonable ordinal notation system for those systems would be best described as only partially solved.

Superficially, Rathjen's $\Psi^\xi_\pi(\alpha)$ (for KP + $\Pi_3$ reflection) looks like C(a,b,c) (or C($\xi,\alpha,\pi$)) in "Degrees of Recursive Inaccessibility". However, while Rathjen does not go beyond ordinary built-from-below construction, he uses a trick that allows $\xi$ to correspond to degree of recursive Mahloness instead of inaccessibility. (Actually, for simplicity his paper uses large cardinals instead of their recursive analogues. Also, while the examples below illustrate his system, I did not verify that the examples hold for all parameters.) Instead of collapsing above c, he collapses below (admissible) $\pi$. In his system, $\pi$ codes both an ordinal above which to do the collapse, as well as the degree to which the resulting ordinal will be a limit of admissibles. If $\pi$ is a successor admissible, the collapse will be an ordinal above $\pi$'s admissible predecessor. If $\pi$ is a successor recursively inaccessible, $\Psi^0_\pi(\alpha)$ will be a limit of admissible ordinals and will be above $\pi$'s recursively inaccessible predecessor, and so on. Admissible ordinals below the least recursively Mahlo ordinal $\Xi(2)$ are obtained using $\Psi^1_{\Xi(2)}(\alpha)$, which works since below $\Xi(2)$ the degree of recursive inaccessibility is in a sense recursive. Now, for another example, suppose that $\pi$ is the least recursively hyper-Mahlo limit of recursively hyper-hyper-Mahlo ordinals above the least recursively hyper-hyper-hyper-Mahlo ordinal $\Xi(5)$. $\Psi^0_\pi(\alpha)$ ranges over non-admissible limits of recursively hyper-hyper-Mahlo ordinals above $\Xi(5)$, $\Psi^1_\pi(\alpha)$ ranges over admissible limits of such ordinals, and $\Psi^2_\pi(\alpha)$ ranges over recursively Mahlo limits of such ordinals.

The downside of Rathjen's approach is the lack of the simple monotonicity that C has, as well additional complexity when trying to extend it (especially beyond $\Pi_3$-reflection) where simple built-from-below does not suffice. Even more recent systems along these lines are fairly complex. For examples, see (Duchhardt 2008) for $\Pi_4$ reflection, (Arai 2015) for $\Pi_n$ reflection, and (Stegert 2010) for stability.

**History of this paper:** In 2005, I discovered the right general form of C, defined a notation system at the level of *a*-recursively inaccessible ordinals, and proposed a schema (a detailed idea) for reaching second order arithmetic (August 7, 2005 version of the paper). In 2006 (January 24, 2006 revision), I defined the stronger ("Degrees of Reflection") notation system, and (June 28, 2009 revision) added the comparison algorithm subsection (the algorithm was previously given implicitly). Finally, in March 16, 2012 revision, I added the main ordinal notation system, intended for second order arithmetic. Afterwards, in 2014 (January 9, 2014 revision), I discovered new structure in the main system and wrote an initial version of a series of updates. In 2015 (April 20, 2015 revision), I proposed a strengthening of the main system. August 24, 2016 revision added "Built-from-below with Passthrough", and around that time, I made many improvements to the paper. October 9, 2016 revision extended the well-foundedness proof for the passthrough system to a stronger system. August 2017 revision added "Using Canonical Definitions" and updated other parts of the paper. November 2018 revision (1) added the correct canonical assignment for Degrees of Reflection (including a correction at the level of $\Pi_n$ reflection), (2) added reflection configurations, and (3) made expository, organizational and other changes.

## 1.3 Goals of Ordinal Analysis

Ordinal analysis of set theories gives us a qualitatively new understanding of the theories and of the infinite ordinals that appear in those theories. A good ordinal notation system captures all ordinals that have a canonical definition in the theory. If a theory is compared to a country, a good ordinal notation system can be viewed as a precise map of the country.

The cumulative hierarchy (or at least, enough aspects of the cumulative hierarchy to acts as a scaffold for large cardinal axioms) and various large cardinal axioms are conceptually simple. We expect they can be captured with reasonably simple combinatorial principles into ordinal notation systems, though we do not yet know whether any of the systems in the paper reach that far.

To state the goals more precisely, for simplicity, let us assume that the theory is naturally consistent with the existence of a definable function f such that $\forall s \exists n \in \mathbb{N}$ s=f(n) (note that, for us, it suffices for s to range over ordinals). If not, we can use a conservative extension of the theory by, for example, a universal Skolem function: $\exists t \varphi(s,t) \Rightarrow \varphi(s,F("\varphi",s))$ ($\varphi$ has two free variables and does not use F) (and an axiom to that effect); alternatively, for many theories, we can add proper classes or their analogues. We are also need a notion of the ordinals; for simplicity, if necessary, let us conservatively extend the theory with ordinals. Sometimes, the most natural extension is nontrivial. For example, for second order arithmetic (or ZF\P) with projective determinacy, we may want to consider Wadge ranks of projective sets instead of just countable ordinals.

In general, it is not clear what canonical means, or even, until we get a suitable ordinal notation system, that it is unambiguous. However, for recursive ordinals, canonical definitions are recursive. If the proof theoretical ordinal of a theory T (T includes $\Pi^1_1$ induction) is sup($a_n$), an example of a non-canonical definition is "sup($a_n$) if each $a_n$ exists, otherwise the largest $a_n$ that exists". Intuitively, unlike provably recursive ordinals in T, this ordinal can be given a more direct definition in an extension of T. Moreover, T is consistent with that ordinal already having a notation $a_n$, albeit for an n that is nonstandard outside of T.

Even for a strong appropriate theory T (extended with Skolem functions if necessary), there is a recursive well-ordering '≺' of $\mathbb{N}$ (with comparison but not well-foundedness provable in T), and an assignment formula φ:$\mathbb{N}$→Ord (inside T, φ may be partial as a function) such that T is consistent with "every ordinal gets assigned a notation by φ, and ordinal comparison agrees with ≺". The reason is that (at least for appropriate T), T is consistent with there being a cut I of $\mathbb{N}$ with every ordinal definable at a definability level in I, with each definability level in a sense transcending previous levels. Now pick a well-ordering '≺' such that T is consistent with '≺' embedding rationals $\mathbb{Q}$ (and hence $\mathbb{Q}^{<\omega}$), and proceed by induction on definability level to assign all ordinals to elements of $\mathbb{Q}^{<I}$. A far-reaching conjecture is that '≺' and φ can be chosen canonically, thus effectively capturing T.

*Notes on canonicity:*
\* It is often easier to skip some canonically definable ordinals, in which case we get only some of the benefit of a full ordinal analysis. A partial remedy is to give canonical assignments for some ordinals (as this paper does). Another partial remedy is to have multiple systems at different strengths each illuminating a different aspect of the theory.
\* Even if a notation system is too weak for a theory, setting gaps "as if it were stronger" can sometimes illuminate aspects of the theory.
\* If we do not build the notation system in L (or a higher canonical inner model), we get an obstacle to canonicity. For example, ZFC does not prove $\omega_1^L < \omega_1$, so an ordinal notation system for ZFC cannot have different terms for $\omega_1^L$ and $\omega_1$. It appears that the resolution is existence of multiple natural assignment formulas φ, but they all agree if V=L (or an analog for stronger theories), and therefore if an ordinal has provably the same definition (for example $\omega_1$) for natural $\varphi_1$ and $\varphi_2$, then it gets assigned the same term by both $\varphi_1$ and $\varphi_2$.

A full completeness (every ordinal has a notation) only occurs in nonstandard models, but perhaps we can have a form of completeness in the true model if we formulate the system as an extensible one. For example, suppose we have a structure with ordinals, and a convention to make every formula φ (alternatively, every φ in a given class) with one-free variable denote an ordinal (for example, let Ord(φ) be the least ordinal satisfying φ, or 0 if there is no such ordinal). Every sound theory T (proving that ordinals are ordered and include 0) gives a c.e. well-founded partial order $[\varphi]_T \le [\psi]_T \Leftrightarrow T \vdash Ord(\varphi) \le Ord(\psi)$. The rank of each formula (in the partial order) is recursive (and worth investigating, but note that some formulas may give definitions that are not canonical for T, and I do not know what the rank of (say) $\omega_1^{CK}$ is for various T), but the canonical assignment for [φ] uses the actual Ord(φ). If $T_1 \subset T_2 \subset T_3 \subset$ ..., converges to the complete true theory (with each $T_i$ c.e.), then at each level, the partial order is c.e., but the limit is the order type of ordinals definable in the structure. An attempt to assign [φ] below

Ord(φ) would cause a problem at some $T_n$. Similarly, if we define an ordinal notation system in a sufficiently extensible manner, the canonical assignment is the lowest appropriately definable assignment that does not break the extensibility. (And even if cannot make the above literally true, it acts as a guiding heuristic.)

*What we want from ordinal analysis of a theory T:*
- A canonical formula φ assigning notations to ordinals.
  -- Ideal: Every canonically definable ordinal in T gets a notation.
  -- Extension: T is consistent with "every ordinal gets assigned a notation by φ". The relative consistency and $\Pi^1_1$-soundness is provable in a weak base theory. T need not prove that if a term is assigned an ordinal, then so do all of its subterms (though it will prove that for every particular numeral).
- A comparison relation '≺' for the system such that provably in a weak base theory (PRA or weaker) '≺' is a polynomial-time computable linear ordering (with a polynomial-time computable domain).
  -- Open ended: proof in a weak base system of various basic properties of '≺'.
  -- Ideal: '≺' is as simple as possible given the strength of T (except that complexity is sometimes justified by convenience, extensibility, or other goals).
- Proof in T that for all m and n such that φ(m) and φ(n) exist, m≺n ⇔ φ(m) < φ(n).
- Proof in a weak base theory that for all n, T proves that φ(n) exists.
  -- Extension (assuming that T is finitely axiomatizable): A proof in T that there is a cut of ℕ such that for k in the cut, φ is total on $I_k$ and exists as a set, where $I_k$ is the set of notations all of whose subterms are below (using '≺') the kth element of a natural fundamental sequence for '≺'. A further extension is to define φ for notations above an arbitrary (possibly within bounds) ordinal, and with T proving that the above holds for every such ordinal (note that the cut of ℕ can be independent of the ordinal).
- Proof that the supremum of notations for recursive ordinals is the proof theoretical ordinal for T. There are several variations on this, but mainly:
  -- Proof in a weak base theory that a $\Pi^1_1$ statement is provable in T iff there is n≺Ω such that the statement is provable in a weak base theory plus transfinite induction on ≺ below n. Here, Ω is the notation for the least nonrecursive ordinal. This should be possible for natural T.
  -- Extension (if the notation system is given above an arbitrary ordinal *a*): Proof in a weak base theory that a $\Pi^1_2$ statement is provable in T iff there is n≺Ω(a) (n uses a placeholder for "a"; Ω(a) is the bound of the recursive part) such that the statement is provable in $RCA_0$ plus well-foundedness of the system up to n for every *a*. This applies to typical T that get their $\Pi^1_2$ strength using constructs like ∀X ∃Y rather than ∃Y. (While $ATR_0$ is needed for basic theory of countable ordinals, a reasonable strong notation system gets $ATR_0$ from $RCA_0$ (Rathjen 2014).)
- (potentially open-ended) A proof theory for T that corresponds to the ordinal notation system. This includes extending the above to higher definability levels.

- Open ended: A detailed correspondence between fragments of T (and their proof theoretical ordinals) and fragments of the notation system. Part of the assignment being canonical is existence of a natural correspondence between n and the proof theoretical ordinal of (essentially) "φ(n) exists". The correspondence is notation system dependent; φ(n) need not be a recursive ordinal, but the proof ordinal is related to n.

We also want the system to be extensible, both generally, and also as a parameterization of the system achieving a form of completeness (see above and also [4.3 Assignment of Degrees](#)) or at least (to a certain extent) capturing the framework that T intuitively provides for its extensions.

The assignment of ordinals can be viewed not just as a formula but also as a story of the ordinals of T. Even without a proof, a detailed correspondence between the notation system and the canonical ordinals (and set-existence principles) of a theory is strong evidence that the system corresponds to the theory.

# 2 A Framework for Ordinal Notations

## 2.1 Definition of the General Notation

In this paper, we define several ordinal notation systems based on the general notion of degree and the corresponding collapsing function C. Key to effective ordinal notation systems is managing complexity and distilling the strength into a simple form. Using C will serve this goal by trivializing syntax (except which forms are standard) and comparison of standard forms while providing a general framework that one can use to achieve arbitrary strength.

**Definition:** Given an arbitrary binary relation on a well-ordered set, which will call a degree (and determine when c has degree $a$), the corresponding *collapsing function C* is defined such that
$C(a, b)$ is the least element above b that has degree $\geq a$. $C(a, b)$ is defined iff such an element exists.
The *strict version of C* requires that $a$ is maximal in $C(a,b)$, with $C(a,b)$ undefined if the least element above b that has a degree $\geq a$ also has a degree $> a$.

*Notes:*
\* C is called a collapsing function because typically, if $a$ is much larger than $b$, then $C(a, b) < a$ and acts as a collapse of $a$ above $b$.
\* The strict version of C is sufficient for standard forms (in the systems we will consider). In this version, existence $C(a,b)$ typically implies that $a$ is in some sense built from below from $<C(a,b)$.
\* The definition also makes sense for sets that are just linearly ordered.
\* A slightly less general presentation (that remains sufficient for us) is that for a well-ordering S, a degree is a partial function $f:S\to S\cup\infty$, and
$C(a,b)=\min(c>b:f(c)\geq a)$ when such c exists. c has degree $a$ iff $f(c)\geq a$. Given a notation system using C (strict version), 0, and Ω (as below), the corresponding choice of f is $f(\Omega)=\infty$, $f(C(a,b)) = a$, and (for definiteness) $f(0)=0$.

Even at this level of generality, C has the following key properties:
* If *a* is maximal in C(a,b) (and in particular for standard C(a,b)), then C(a,b) <C(c,d) iff C(a,b)≤d ∨ (b<C(c,d) ∧ a<c).
* The order type of terms built using 0 and C is bounded by $\varepsilon_0$. Proof sketch: The canonical choice of C here is $C(a,b)=b+\omega^a$, and every other choice is essentially a quotient of this choice.
* Consider an arbitrary notation system built using 0 (the least ordinal), Ω (a large ordinal), and C (strict version) with b<Ω ⇒ C(a,b)<Ω.
  - The order type of terms above Ω is bounded by $\varepsilon_{\Omega+1}$.
  - The order type of terms below Ω can be
    -- an arbitrarily large recursive ordinal for a recursive system
    -- an arbitrarily large countable ordinal for a nonrecursive system
    -- an arbitrarily large ordinal $<\kappa^+$ for a nonrecursive system if all ordinals ≤κ (κ<Ω) are added as constants and C(a,b)>κ.
  - The comparison of valid terms (terms that we are promised to be valid) is recursive (using the formula above) and independent of the notation system. For the extension above κ, recursive comparison holds given comparison of the constants <κ used in both terms.
  - Furthermore, assuming the system always uses C(a, c) in place of C(a,C(b,c)) if b<a (note that these terms are always equal for strict C), then comparison of standard terms is lexicographical if the terms are written in postfix form, with 'C'<0<Ω. For example, C(C(0,Ω),0) would be written as 0Ω0CC, and (assuming standardness) 0ΩC < 0Ω0CC < 0Ω0CC0C < Ω. This also applies to the extension above κ if we use C(a,κ) in place of C(a,b) for b<κ.

To ensure that enough terms are standard (without restricting how strong the system can be), we will only use degrees that are well-behaved as follows.

**Definition:** A *degree* for a well-ordered set S is a binary relation on S such that
  1. Every element c∈S has degree $0_S$ (the least element of S). $0_S$ only has degree $0_S$.
  2. For a limit *a*, c has degree *a* iff it has every degree less than *a*.
  3. For a successor *a'=a+1*, either of the following holds:
      ○ An element has degree *a'* iff it is a limit of elements of degree *a*.
      ○ There is a limit element d ≤ *a* such that for every c in S, c has degree *a'* iff it has degree *a* and either c ≤ d or c is a limit of elements of degree *a* (or both).

*Note:* The third condition can be equivalently written as ∀a ($C_{a+1}$ = lim($C_a$) ∨ ∃d ∈ lim(S)∩(a+1) $C_{a+1}$ = lim($C_a$) ∪ ($C_a$∩(d+1))), where S is identified with an ordinal (so a+1 consists of ordinals ≤a), $C_a$ is the set of elements that have degree *a*, and lim is limit points.

Intuitively, a degree indicates how closed the ordinal is. Degrees are fine-grained in that degree a+1 includes all limit points of the set of ordinals of degree *a*. For typical notions of degrees, the intuitive meaning of c having degree *a* depends not only on how large *a* is but also on the nature of the gap between c and *a*; this is necessary to preserve monotonicity. Given a typical ordinal notation system, the representation of *a* in terms of ordinals <c can be viewed as a reflection configuration with ordinals <c being parameters. Intuitively, if we increase c past certain subterms of *a*, the meaning of these subterms changes from being building blocks of the reflection configuration to that of parameters. Also, not all representations are valid as reflection configurations. Some increases in *a* correspond to diagonalization: We may have an increasing sequence of degrees $a_1, a_2, ..., a_{diag}$ with c having degree $a_{diag}$ meaning c (if it is not too large relative to $a_{diag}$) having degree $a_c$. To preserve monotonicity; $a_i$ may only be used when i<c. For i≥c, we skip $a_i$ by making it equivalent to $a_{diag}$ (and hence to $a_i+1$ and higher degrees <$a_{diag}$) for c.

Besides C, we might only need 0 and Ω (as above). More generally, however, a *partial ordinal notation system* is a partial mapping O of ordinals into finite sequences of symbols and nonrepresentable ordinals (ordinals not in the domain of O).

**Definition:** For a partial ordinal notation system O and a collapsing function C corresponding to a degree, the standard combination of the notations is defined as follows:
- If an ordinal *a* is in the domain of O, then represent *a* using O and the representations of ordinals in O(*a*).
- Otherwise if *a* = C(b, c) where b is maximal and c is minimal (that is b'>b ⇒ C(b', c) > C(b, c), and c' < c ⇒ C(b, c') < C(b, c)), then represent *a* through C(b, c).
  Exception: If O includes all ordinals ≤κ as constants, it is often convenient to require that in C(a,b), b≥κ. This helps with embedding notation systems built above different κ.

In general, we can stack ordinal notation systems on top of each other. For example, we can use one system for ordinals below Ω, and another system for representing ordinals above Ω in terms of ordinals below Ω. A well-behaved combined ordinal notation system should satisfy the following (this is not a definition and can be skipped by the reader; the conditions are most relevant when an ordinal notation system is stacked on top of another):
1. Let Ω be the least ordinal that is above all ordinals that are represented by C-terms (a C-term is just C(…, …)). C(a, b) is defined whenever *a* and b are representable and min(a, b) < Ω (except possibly when b≥Ω and b+$ω^a$ is above all ordinals in the notation system).
2. Let *a* be the least ordinal represented by a C-term and Ω as above. If a ≤ C(b,c) < Ω and *b* is maximal and c is minimal in C(b, c), then C(b, c) is represented by C-term (C(b,c)) except possibly when c=0 and b is not represented by a C-term.
3. The standard representation of an ordinal is the unique representation using (C and 0) with the least number of C-terms (except that we can have C(a,0)=C(a,b) even if b is not a C-term). This should hold even if all representable ordinals below a particular one were to be added to O as constant symbols.

An example is the notation system below $\varepsilon_0$ generated using 0 as the least ordinal and $C(a, b) = b+\omega^a$. The standard representation is analogous to Cantor normal form. Another example is generated by 0, $x\to\varepsilon_x$, and C with $C(a, b) = b+\omega^a$.

The notation systems in this paper can also be built on top of an arbitrary ordinal by representing all ordinals below that one by a constant.

More precisely/generally, a notation system above a generic ordinal specifies standardness and comparison given syntactic form and comparison for inputs. For example, pairing gives a notation system for $\Omega^2$ given a system for $\Omega$, and comparison of (a,b) with (c,d) depends only on how a,b,c,d compare with each other. (Note: While not necessary here, it may be convenient to relax the definition of ordinal notation system above a generic ordinal by allowing comparison with 0 and possibly other queries.) Notation systems in this paper can generally be defined above a generic ordinal. Furthermore, the general form of some constructs is a mapping that given an ordinal notation system above a generic ordinal, turns it into a stronger one. A key value of defining the notation system above a generic ordinal is that it can be used to characterize not just $\Pi^1_1$ but $\Pi^1_2$ strength of typical theories.

## 2.2 Basic Properties

Using the general definition, we can prove a variety of properties of C. Without loss of generality, we will assume in the following that the well-ordered set S is an ordinal.
If b < a, then having degree *a* implies having degree *b*.
If $0 < a \leq b$, then *b* has degree *a* iff there are *c* and *d* with $d \geq a$ and $b = c+\omega^d$.
$C(a, b) = b+\omega^a$ iff $C(a, b) \geq a$.
$C(a, b)$ is monotonic in *a* and *b* and is continuous in *a*.
$C(a, b) > b$.
If $b = C(c, d)$ where *c* is maximal, then b is minimal in $C(a, b)$ iff a≤c.
In general, if *e* is represented in the standard form as $C(a_1, C(a_2,... C(a_m, b)...))$ and f as $C(c_1, C(c_2,... C(c_n, b)...))$, then $e < f \Leftrightarrow (a_m, ... a_2, a_1) < (c_n, ..., c_2, c_1)$ with the comparison done in lexicographical order. Also, $a_m \geq ... a_2 \geq a_1$.
If *a* is maximal in $C(a, b)$ and $b \leq d$, then *a* is maximal in $C(a, d)$.

*Comparison*:
If *a* is maximal in $C(a, b)$, then $C(a, b) < C(c, d)$ iff $C(a, b) \leq d$ or ($b < C(c, d)$ and a < c).
Since the standard form requires *a* to be maximal in $C(a, b)$, the above equivalence leads to a polynomial time comparison algorithm for ordinals in the notation system assuming that we can compare *a* and *b* when "a" is a constant or is otherwise not a C-term.
Note that to evaluate $C(a, b) < C(c, d)$, *c* need not be maximal in $C(c, d)$. In fact, recursive definitions of maximality (and standard form) will often involve such comparisons.

## 2.3 Bachmann-Howard Ordinal

We illustrate the framework by defining a notation system at the level of the Bachmann-Howard ordinal, which is the proof ordinal of KP. The notation system will consist of C (as above) and two constants: 0 (the least ordinal) and $\Omega$ with C(a, b) < $\Omega$ if b < $\Omega$. At this point, the notation system is fully specified except for the definition of when *a* is maximal in C(a, b). On the other hand, at this point the system is fully general in that by choosing the maximality condition, we can reach arbitrarily high recursive ordinals. Note that the definition of maximality cannot affect comparison of ordinals in the standard representation, but only determines which representations are standard. To reach the Bachmann-Howard ordinal, we use the built-from-below condition: define *a* to be maximal in C(a, b) iff its standard representation (equivalently, some representation) in terms of ordinals ≤b (equivalently, ordinals <C(a, b)) only uses ordinals below *a*.

*Intuitive Explanation:* The built-from-below condition intuitively captures the idea of $\Omega$ being admissible. A subsystem of the notation system obtained by setting an upper bound on subterms used corresponds to a restriction on the theory and thus to a provably recursive well-ordering. Conversely, the theory corresponds to the union of its restrictions, and the limit/sum of these subsystems corresponds to the proof ordinal (which in a sense, corresponds to $\omega_1^{CK}$), and this limit is precisely the order type of built-from-below ordinals. This way, the ability to collapse built-from-below ordinals into ordinals below $\Omega$ corresponds to recursive ordinals being below $\Omega$. The above also generalizes to terms built above an ordinal (and for extensions of the notation system, to successor admissible ordinals). Also, ordinal exponentiation above $\Omega$ corresponds to quantifiers in induction instances, analogously to PA and $\varepsilon_0$.

For C($\Omega$*a+b, c) (a, b, c < $\Omega$), either $\Omega$*a+b or $\Omega$*(a+1) or $\Omega^2$ is maximal. C($\Omega^2$, c) is the least ordinal in the range of $\Gamma$ above c. If $\Omega$*a+b is maximal and a > 0, then C($\Omega$*a+b, c) is point number $\omega^b$ above c in the range if the *a*th fixpoint function ($\varepsilon$ is the first fixpoint function).

Let H(a, b) is defined as the least set of ordinals including $\Omega$ and all ordinals below *b* and closed under x,y → C(x, y) where x<a, then
For b < $\Omega$, C(a, b) = H(a, b+1) ∩ $\Omega$ = min(c>b: c = H(a,c) ∩ $\Omega$). (If b ≥ $\Omega$, then C(a, b) = b+$\omega^a$.)
Note: H would be unchanged if "closed under x,y → C(x, y)" is replaced by "closed under x,y → C(x, y) where x is maximal and y is minimal in C(x, y)".

To convert C(a, b) to the standard form, convert *a* and *b* to the standard form, maximize *a* and minimize *b*. b can be minimized using C(x, C(y, z)) = C(x, z) if y < x and y is maximal in C(y, z) (and if b is $\Omega$ or 0, then it is already minimal). To maximize *a*, if a < $\Omega$ and not already maximal, then just replace *a* with $\Omega$. If a ≥ $\Omega$ (and not already maximal), then look at the standard representation of *a* in terms of ordinals below $\Omega$. Find the rightmost ordinal instance whose standard representation in terms of ordinals below b+1 uses an ordinal above *a*. Replace that instance with $\Omega$ and delete all ordinals (including $\Omega$) to the left of it, removing extraneous Cs (by recursively replacing C(,x) with x and deleting C(,)).
For example, if b = 0 and a = $\Omega*\omega+\Gamma_0*2+1$ = C(0, C($\Gamma_0$, C($\Gamma_0$, C(C(0, $\Omega$), $\Omega$)))), then out of (0, $\Gamma_0$, $\Gamma_0$, $\Gamma_0$, 0), delete the first two and convert the third one to $\Omega$: a $\to$ C(, C(, C($\Omega$, C(C(0, $\Omega$), $\Omega$)))) $\to$ C($\Omega$, C(C(0, $\Omega$), $\Omega$)) = $\Omega*(\omega+1)$.

## 2.4 One Variable C

To allow more readable ordinal representations, given a definition of C, let us define $C_1$ as
$C_1(a) = C_1(a,0)$
$C_1(a,b) = C(a_1,C(a_2,... C(a_n,b)..))$ for the least ordinal $\omega^{a_n}+...+\omega^{a_2}+\omega^{a_1} \geq a$ such that each $a_i$ is maximal in $C(a_i,... C(a_n,b)..)$.
(Note that the definition is unchanged by requiring $a_n \geq ... \geq a_2 \geq a_1$.)
Equivalent definition of $C_1$:
$C_1(a) = C_1(a,0)$
$C_1(0,b) = b$
if $d < \omega^{a+1}$ and $a$ is maximal in C(a,b): $C_1(\omega^a+d,b) = C_1(d, C(a,b))$
if $d < \omega^{a+1}$ and $a$ is not maximal in C(a,b): $C_1(\omega^a+d,b) = C(a,b)$
**Notes:**
* The Python module uses "C" in place of 1-variable $C_1$, C(a, base=b) in place of $C_1$(a, b), and "C2" in place of the 2-variable C.
* If one needs two-argument $C_1$, it may be helpful to switch notation so that comparison remains left-to-right. For standard representations, the base b is the more significant term than a.

Thus,

- $C(a_1,C(a_2,... C(a_n,b)..))$ (if each $a_i$ is maximal) = $C_1(\omega^{a_n}+...+\omega^{a_2}+\omega^{a_1},b)$
- $C(a, b) = C_1(\omega^a, b)$.
- If b is fixed, $C_1$ is continuous and monotonic and onto ordinals $\geq$b. Also, $C_1$ is monotonic in b.
- If $a$ is not above the least fixpoint of exponentiation base $\omega$ above b, then $C_1(a,b) = b+a$. In all cases, $C_1(a,b) \leq b+a$
- If $a_1 < a_2$ and $a_1+b$ is maximal in $C_1(a_1+b,c)$ and $a_2$ is maximal in $C_1(a_2,c)$, then $a_2+b$ is maximal $C_1(a_2+b,c)$.
- $C_1(a,b) = C_1(C_1\text{inv}(b)+a)$ where $C_1\text{inv}(b)$ is the largest x such that $C_1(x)=b$ provided that for $\exists x\ b=C_1(x)$.

Formal definitions of ordinal notation systems tend to be simpler using C than $C_1$ because C already takes care of addition and exponentiation, while with $C_1$, addition and exponentiation would have to be defined separately. However, $C_1$ allows much more readable representations because a single $C_1$ can combine a long chain of Cs, and because key properties $C_1(a)$ can be read off CNF (for appropriate base) representation of *a*. Also, combining $C_1$ with extended CNF (always using the maximum possible fix-point x→$\omega^x$ as the exponent base) has the benefit that the most significant terms are always to the left, which allows easy comparison.

For example (with C as defined in the Bachmann-Howard Ordinal subsection above):

$\varepsilon_{\omega^2+\omega+3} = \varphi(1, \omega^2+\omega+3) = C_1(\Omega*(\omega^2+\omega+3)) = C(\Omega,C(\Omega,C(\Omega,C(\Omega+1,C(\Omega+2,0)))))$

$\varphi(1,\varphi(\varphi(1, \omega^2+\omega+3),5)+6) = C_1(\Omega^{C_1(\Omega*(\omega^2+\omega+3))*4+\Omega*6})$

As seen above, another key simplification for ordinal representations is that typically $C_1$ can be just a one-argument function while two arguments are required for C.

Also, $C_1$ allows using integers (including number of repetitions with the same *a*) in binary instead of unary; this shortcut can also be obtained using (for example) $C_n(i,a,b) = C(a,...,C(a,b)..)$ (i *a*s; 'n' is a character and not a variable).

The conversion between using C and using $C_1$ is straightforward. Typically, in the standard form (and assuming that $C_1(a)$ is representable), "$C_1(a)$" is used iff *a* is maximal and there is no other representation of $C_1(a)$; in particular, if "$C_1(a)$" is standard, then $C_1(a) = \omega^{C_1(a)}$. Typically, two-argument $C_1$ is not used in standard representations. Comparison of standard forms is based on monotonicity of $C_1$ and on $C_1(a)$ being a fix-point of exponentiation base ω. For typical $C_1$, assuming that *a* is maximal in $C_1(a)$, $C_1(a+\omega^b) = \omega^{C_1(a+\omega^b)} \Leftrightarrow \exists x \le b\ x=\omega^x > C_1(a)$.

## 2.5  Reflection Configurations

While the use of degrees allows a simple formalization, it obscures a fundamental symmetry in our notation systems. Namely, every pair (b,a) of representable ordinals defines a special function r(a,b) (with r(a,b)(b) = a) describing how big is *a* relative to b, which can then be applied to other appropriate ordinals.

**Definition:** The *reflection configuration* of a term *a* above b is obtained by
- taking the standard representation of *a* in terms of ordinals/terms ≤b (so subterms ≤b are not parsed for subterms C(c,d)>b),
- in every subterm C(c,d)>b with d<b, replacing d with b
- replacing b with x, and adding λx to the front of the configuration. The variable name x is immaterial as long as x is not used otherwise.
The *constants* used are the outermost subterms <b in the configuration (i.e. a constant (<b) is not parsed for further constants).
The *upper bound* will typically be Ω or infinity, but formally the lesser of the least constant term >b used in the notation system and the term (if any) above all ordinary C-terms.
*Note: a* and b can be ordinals or terms in the notation system used. We are not restricting *a* and b provided that *a* has a standard representation in terms of ordinals ≤b.

The notation systems we will define will have the property that if A is a reflection configuration (of *a* above b), and d is an ordinal above the constants used and below the upper bound, then there is unique c such that A is the reflection configuration of c above d. In other words, given *a* and b and d, c is to d what *a* is to b. Furthermore, c is strictly monotonic in *a* and monotonic in d. In general, c will depend on the notation system used. However, if b and d have sufficiently strong reflection properties relative to the notation systems, c will be independent of the notation system, so the differences in notation systems are in how to handle reflection configurations above ordinals without sufficient reflection properties.

A single *a* has a finite set of reflection configurations, with configuration transitions (of the configuration of *a* above b) happening at subterms of *a*. The final reflection configuration of *a* will just use *a* as a constant (i.e. λx.a). In terms of the structure of the systems, reflection configurations can be said to be more important than degrees.

The comparison of configurations is simply the standard comparison with x treated as an ordinal above the constants but below the upper bound. In our systems using C,0,Ω, the lexicographical comparison (of postfix forms) works provided that we are given whether x<Ω (and similarly with other special constants for systems that use them). We will typically only care about configurations for x<Ω, but one can specify it explicitly using λx<Ω. Note that λx.x is the least nonconstant configuration.

*Example:* For the system for the Bachmann-Howard ordinal, the configuration for $\varepsilon_0$ above any lower ordinal is λx.C(Ω,x), and in line with this value, the least point in the range of ε above b<Ω equals λx.C(Ω,x)(b) = C(Ω,b).

The reflection configuration of *a* above <b is the reflection configuration of *a* above a sufficiently large term <b, except that for nonlimit b, we treat it as the configuration above an imaginary ordinal above all of b's predecessors. Thus, the only difference with using '<' is when *b* is a subterm of *a*.

Of particular importance is the reflection configuration of *a* above C(a,b) (standard term). If the configuration is λx.A(x), the configuration of C(a,b) above <C(a,b) is λx.C(A(x),x), which characterizes how big/reflective C(a,b) is relative to lower ordinals. Note that both configurations are also valid above every x<C(a,b) that is above all outermost subterms of C(a,b) that are <C(a,b).

*Note on testing for standard forms:* If a condition for standardness relies on the configuration of *a* above C(a,b), then in testing it, use standard comparison to determine which terms are <C(a,b). Then apply the condition to check for standardness. For example, C(C(Ω,0),0) is only standard if C(C(Ω,0),0) < C(Ω,0), so the configuration to check (for a=C(Ω,0)) is λx.C(Ω,x) rather than λx.C(Ω,0).

**Notations related to reflection configurations**

To avoid repeating ourselves, here are some notations that make defining particular systems easier. The reader is free to skip this part and come back when encountering a definition that uses reflection configurations.

In testing C(a,b) for standardness with *a* and b standard and b minimal, the remaining condition is about *a*.

* $T_a$ is the subterm/parse tree of 'a': $T_a$ is the set of subterms of 'a', and for x and y in $T_a$, x⊏y means that x is a proper subterm of y. Identical terms at different positions are distinguished (in the quantifiers and in ⊏).
*Note:* $T_a$ will be useful in stating the conditions in first order logic.
* d=parent(c) means c⊏d ∧ ¬∃e c⊏e⊏d, but as a special case, parent(a) is C(a,b).
* d=$parent_{<Ω}$(c) means c⊏d ∧ d<Ω ∧ ∀e (c⊏e⊏d) e≥Ω, or if there is no such d, $parent_{<Ω}$(c) is C(a,b). In other words, if ordinals ≥Ω are merely syntactic constructs, d is the real parent of c.
* r(c,d) is the reflection configuration of c above d.
* r(c) is r(c,parent(c)), annotated with the position of c inside Ω; the standard comparison is used for C(a,b) even if it is nonstandard.
* $r_{<Ω}$(c) is as above using r(c,$parent_{<Ω}$(c)).
*Note:* r(c) or (as appropriate) $r_{<Ω}$(c) are the relevant configurations for extending built-from-below to reflection configurations.
* d∈r(c) iff d⊑c ∧ ∀z (d⊏z⊑c) z≥parent(c) (and similarly with $r_{<Ω}$(c); the point is that r(c) uses d, possibly as a leaf term)
* r(e)⊆r(c) iff ∀z (z∈r(e)⇒z∈r(c))

For a reflection configuration *a*, the subterm tree $T_a$ is defined as if it were an ordinal, except that the constants <x are not parsed further. r is defined as above (without fallback to C(a,b)) by replacing x with a fictious (sufficiently large below the upper bound) ordinal $x_0$, computing the configuration, and replacing λx with λ$x_0$ λx. If a configuration already uses λ$x_0$, then use $x_1$ (like $x_0$), and replace λ$x_0$ λx with λ$x_0$ λ$x_1$ λx; and so on. The comparison remains standard (lexicographical in postfix form with (if x<Ω) C < 0 < $x_0$ < $x_1$ < $x_2$ < ... < x < Ω). For x<Ω, $r_{<Ω}$ is defined analogously.

Comparison between configurations for x<Ω and x≥Ω is by default undefined. However, in our systems using C,0,Ω it will be convenient to compare them as follows: For the configuration for x<Ω, replace Ω with Ω' (except in the constants below x) and then use the usual (lexicographical for postfix forms) comparison with C < 0 < Ω < x < Ω' (actually, any use of Ω' will make the configuration greater). For example, λx<Ω.5 = λx≥Ω.5 < λx>Ω2.Ω2 < λx<Ω.x < λx<Ω.Ω. See [4.4 A Step towards Second Order Arithmetic](#) for an extenson of this for certain systems.

# 3    Degrees of Recursive Inaccessibility

In this section, we define an ordinal notation system for rudimentary set theory plus "for every ordinal *a*, there is recursively *a*-inaccessible ordinal". Instead of rudimentary set theory, we can use a stronger theory like $\Pi^1_1$-$CA_0$.

*Note on rudimentary set theory:* Rudimentary set theory consists of extensionality, foundation, empty set, and closure under rudimentary functions (as a schema, or using a selected set of rudimentary functions). Assuming infinity, it is conservative over $ACA_0$, but it is more general in that it does not assume countability, and (in the other direction) handles sets without assuming $ATR_0$. Every level of the Jensen hierarchy for L satisfies rudimentary set theory, and conversely it can be used to define L and be the basis of a natural axiomatization of $J_a$. Outside of L, the theory is too incomplete to fully function as set theory, but one possible remedy (assuming choice but remaining conservative over $ACA_0$) is to require that every set S is in L[X] for some X⊂$\omega_\alpha$ where X depends on S, and $\omega_\alpha$ is the cardinality in the transitive closure of S (assuming the closure is infinite).

## 3.1 Definition of the Notation

The notation is built from the constant 0 (the least ordinal) and a function C as follows:

1. C(a, b, c) is the least ordinal *e* of admissibility degree *a* that is above *c* and is not in H(b, e).
2. H(b, e) is the least set of ordinals that contains all members of *e*, and is closed under h, i, j → C(h, i, j) where i < b.
3. If an ordinal *e* is of admissibility degree *a*, then C(h, i, j) < e whenever h < a and j < e. 0 is of admissibility degree 0.

A canonical notion of admissibility degrees is:
3'. Ordinals of admissibility degree a+1 are the recursively *a*-inaccessible ordinals and their limits.
For limit *a*, having admissibility degree *a* is the same as having every admissibility degree below *a*.
(An ordinal is recursively *a*-inaccessible iff it is admissible >ω and for every *b* < *a* is a limit of recursively *b*-inaccessible ordinals.).

1, 2, and 3 define the comparison relation for ordinals in the notation. 1, 2, and 3' uniquely fix C and imply 3.

**Interpretation:** C(a, b, c) is the least ordinal above *c* of admissibility degree *a* that is not reachable from below using "collapses" of ordinals less than *b*. If C(a, b, c) < b, then C(a, b, c) may be viewed as a collapse of *b* above *c*.
**Note:** The pair (a, b) may be viewed as a degree in a similar sense as degrees in the Definition of the General Notation subsection.

## 3.2 Comparison Relation

A polynomial time comparison algorithm is obtained as follows:

1. In C(a, b, c), treat (a, b) like an ordinal and use the comparison relation in the Definition of the General Notation subsection. Specifically, if *a* and *b* are maximal in C(a, b, c), then C(a, b, c) < C(d, e, f) iff C(a, b, c) ≤ f or c < C(d, e, f) ∧ (a, b) < (c, d).
   Note: (a, b) < (c, d) iff a < c, or a = c and b < d.
2. *a* is always maximal in C(a, b, c).
3. *b* is maximal in e = C(a, b, c) iff b is in H(b, e).
4. *k* is in H(b, e) iff k < e, or k = C(g, h, i), *h* is maximal, *i* is minimal, each of g, h, i is in H(b, e), and h < b.

Instead of 3 and 4, we can use the following: b is maximal in C(a, b, c) iff it has a representation (equivalently, standard representation) in terms of ordinals below C(a, b, c) (used as constants) using only ordinals below b.

### 3.3 Examples

Here are some examples where C(a, b)=C(0, a, b), and $x^+$=C(1, 0, x) (the next admissible ordinal above x):
C(0, x) = x+1,
C(b, c) = H(b, c+1) ∩ $c^+$, so
C(b, c) is the least ordinal not in H(b, c+1).
However, C(1, $0^+$, 0) is greater than C(1, $\varepsilon_0$, 0), the least ordinal of admissibility degree 1 that is not in H($0^+$, 1).
C(b, c) = min(c+$\omega^b$, C($c^+$, c)) if b < $c^+$,
C($a^+$, a) is the least fixed point of x → $\omega^x$ above *a*.
{x: c < x < C(a, b, c) and x has admissibility degree *a*} has order type 0 if b=0 and min($\omega^b$, C(a, b, c)) otherwise.
C(a, C(a+1, 0, b), b) is the least ordinal c above b such that there are c ordinals <c of admissibility degree *a*.
C($(0^+)^2$, 0) = $\Gamma_0$.
C($(0^+)^\omega$, 0) is the small Veblen Ordinal, and C($(0^+)^{0^+}$, 0) is the large Veblen ordinal.
C($0^{++}$, 0) is the Bachmann-Howard ordinal.
C(C(1, 1, 0), 0) is the proof-theoretical ordinal of $\Pi^1_1$-$CA_0$.
C(C(1, 1, 0)$^+$, 0) is the ordinal for $\Pi^1_1$-CA + TI.
C(C(2, 0, 0), 0) is the ordinal for $\Pi^1_1$ Transfinite Recursion (with induction limited to sets).
C(C(2, 0, 0)$^+$, 0) is the ordinal for KPi.
C(C(3, 0, 0)$^+$, 0) is the ordinal for KP + a proper class of recursively inaccessible cardinals.

## 4 Degrees of Reflection

### 4.1 Definition of Degrees of Reflection

In this section, we present an ordinal notation system for KP + $\{\Pi_n \text{ reflection}\}_n$, which should have the same strength as $ACA_0$ + lightface $\Pi^1_2$ comprehension. The version using Veblen Normal Form reaches beyond existence of $\kappa^+$-stable $\kappa$.

**Definition:**
**Syntax:** 0 (the least ordinal), C (binary function), Ω (large ordinal), **O**.
**Condition for O: O** is a notation system for ordinals >Ω in terms of ordinals <Ω. Ordinals <Ω are treated as given in **O**: Comparison and validity of **O**-terms depends only on the term structure and on comparison of the ordinals <Ω used by the terms.
**Definition of C:** The notation uses C(a, b) for b < Ω, which for ordinals in the notation implies C(a, b) < Ω, and **O** for larger ordinals, and 0 for the least ordinal. C is as given in "A Framework for Ordinal Notations". This fully specifies the system (if we do not use gaps) apart from maximality of *a*. *a* is maximal in C(a, b) iff for every ordinal d in the **O** representation of *a* (d<Ω), the following holds. The standard representation of d does not use ordinals that are below Ω but greater than d, excluding instances in the scope of an ordinal less than C(a, b). (If d is ...f..., and f < C(a, b), then do not parse f for ordinals larger than d. If a < Ω, then d is a.)
**Canonical choice of O:** For definiteness, let **O** be Cantor Normal Form base Ω. Other choices include C with C(a, b) = b + $\omega^a$ (b ≥ Ω; in the standard form b = Ω or b ≥ max($\omega^a$, Ω)), or for a stronger system, Veblen Normal Form.

**Notes:**
* For a more explicit phrasing, see Comparison Algorithm below.
* When C is used for **O**, the system is superficially similar to (but weaker than) the n=2 notation system in [5 Main Ordinal Notation System](5 Main Ordinal Notation System).
* A natural strengthening is [8.2 Degrees of Reflection with Passthrough](8.2 Degrees of Reflection with Passthrough).
* Compared to using C, using CNF base Ω has the same strength but allows additional terms. For example, if a = Ω+d with d = $d_1$+...+$d_n$ <Ω (each $d_i$ is additively indecomposable and ($d_i$:1≤i≤n) is nonincreasing), we only check d rather than checking the representation of each $d_i$ for ordinals >$d_i$. However, if the representation of *a* using CNF base Ω does not use ordinals >C(a,b) that are not fix-points of x→$\omega^x$, then C(a,b) is the same in both systems (but converting C(a,b) still requires conversion of *a* and b).

**Comparison Algorithm**

Here, we explicitly state the comparison algorithm and the conversion to the standard form (when C is used for **O**).

**Syntax:** Two constants (0, Ω) and a binary function C.
**Comparison:** For ordinals in the standard representation written in the postfix form, the comparison is done in the lexicographical order where 'C' < '0' < 'Ω': For example, C(C(0,0),0) < C(Ω,0) because 000CC < 0ΩC.
**Standard Form:**
0, Ω are standard

"C(a, b)" is standard iff
1. "a" and "b" are standard,
2. b is 0 or $\Omega$ or "C(c, d)" with a≤c, and
3. $\forall x \in T_a \, \forall y \sqsubset x \, (x<y<\Omega) \, \exists z \sqsupseteq y \, (z<\Omega) \, (z \sqsupseteq x \lor z < C(a, b))$.

Here $T_a$ is the parse tree of 'a': $T_a$ is the set of subterms of 'a', and for x and y in $T_a$, x⊏y means that x is a proper subterm of y; identical terms at different positions are distinguished (in the quantifiers and in ⊏).
Note that $\exists z \sqsupseteq y \, (z<\Omega) \, z \sqsupseteq x$ simply means that x is not in the representation of *a* in terms of ordinals <$\Omega$ (with z⊒y redundant).
**Note:** The '<' uses the above-described lexicographical order (including for z < C(a, b)).

*When using an arbitrary **O** for ordinals >$\Omega$:* Condition (3) applies verbatim. (The term structure above $\Omega$ is irrelevant as long as the ordinals used <$\Omega$ are unaffected.) The other conditions are standard for C and can be handled accordingly. (In condition (2), exclude b=$\Omega$. Use standard comparison for C. Comparison and standardness of **O**-terms is per the definition of **O**.)

*Conversion of nonstandard forms:* To convert 'C(a, b)' to the standard form (when C is used for **O**), first convert 'a' and 'b'. Next, recursively minimize b by replacing it with d for as long as b is C(c, d) and c<a. If we are not done, then perform a right-to-left (with functions written in prefix form) preorder traversal of $T_a$ until we find the first *y* that violates the above condition. Replace *y* with $\Omega$, and delete everything in 'a' to the left of 'y' (except for the required number of 'C'). Convert 'a' to the standard form. 'C(a, b)' is now a standard representation.

## 4.2 Examples and Additional Properties

**Examples**
Here are examples using a canonical **O**, and setting gaps in the canonical way.
$C(\Omega*a+b, c) = C(a, b, c)$ of the previous notation (for a, b, c representable in that notation; **O** uses CNF base $\Omega$).
Ordinals (below $\Omega$ and with a<$\Omega$) of degree $\Omega^{a+1}$ are recursively *a*-Mahlo ordinals (and their limits), except that if $c=C(\Omega^{a+1},b) \leq a$, c is recursively c-Mahlo.
$C(C(\Omega, C(\Omega^2, 0)), 0)$ is the ordinal for KPM (KP + the universe is recursively Mahlo).
In general, for many appropriate conditions F (including $\Pi_n$ reflection), the ordinal for KP + "the universe is F" is $C(C(\Omega, a), 0)$ (which equals $C(\varepsilon_{a+1}, 0)$) where *a* is the least ordinal such $L_a$ satisfies or can be forced to satisfy F.
However, the ordinal for KP + "for every ordinal *a*, there is recursively *a*-inaccessible ordinal" is $C(C(\Omega^2, 0) + \varepsilon_{a+1}, 0)$ since $C(\Omega^2, 0)$ is needed to reach a $= C(\Omega*C(\Omega^2, 0) + \Omega, 0)$, the least ordinal *a* that is recursively *a*-inaccessible.
Ordinals below $\Omega$ of degree $\Omega^\Omega$ are $\Pi_3$ reflecting ordinals (and their limits).
Ordinals below $\Omega$ of degree $\Omega^{\Omega^\Omega}$ are $\Pi_4$ reflecting ordinals (and their limits), and similarly with $\Pi_n$ reflection.

An ordinal that is a limit of ordinals of maximum degree *a* has degree a+b iff it is a "level" b limit of ordinals of degree *a*. For example, $C(\Omega^\Omega+\Omega^2+\Omega+1, 0)$ is the least limit of admissible limits of recursively Mahlo limits of $\Pi_3$ reflecting ordinals.
An analogous property holds for ordinals of degree a*b and $\Pi_2$ reflection. For example, $C(\Omega^{\Omega^\Omega+\Omega+1},0)$ is the least ordinal that is $\Pi_2$ reflecting onto $\Pi_3$ reflecting ordinals that are $\Pi_2$ reflecting onto $\Pi_4$ reflecting ordinals.
Similarly, $C(\Omega^{\Omega^{\Omega^\Omega+\Omega+1}},0)$ is the least ordinal that is $\Pi_3$ reflecting onto $\Pi_4$ reflecting ordinals that are $\Pi_3$ reflecting onto $\Pi_5$ reflecting ordinals, and similarly for $\Pi_n$ reflection.

*Levels of Stability (using Veblen normal form for **O**)*:
$C(\Omega^{\varepsilon_{\Omega+1}+\Omega^\Omega+1},0)$ is the least ordinal that is $\Pi_2$ reflecting onto $\Pi_4$ reflecting ordinals that for every finite n are $\Pi_2$ reflecting onto $\Pi_n$ reflecting ordinals, and similarly with other expressions below $C(\varepsilon_{\Omega+2},0)$.
$C(\varepsilon_{\Omega*2},0)$ is the least 1-stable ordinal.
$C(\varphi_\Omega(1),0)$ is likely the least ordinal *a* that is $a^+$-stable (equivalently, the least $\Pi^1_1$-reflecting *a*).
$C(\varphi_{\Omega*\Omega+\Omega}(0),0)$ is essentially the least ordinal *a* corresponding to iterating $\Pi^1_1$ reflection $a^{+}*a^{+}+a^{+}$-times, and analogously for other ordinals.
**Notes:**
* $a^+$ is $C(\Omega, a)$ (a<$\Omega$); *a* is *b*-stable iff $L_a \preceq_{\Sigma_1} L_{a+b}$; *a* is stable iff it is *b*-stable for every *b*.
* If one extends **O** to stronger recursive systems (as opposed to systems imitating nonrecursive structure) while wishing for $\Omega$ to remain the same, a natural choice for $\Omega$ is the least non-Gandy ordinal, equivalently the least $\Sigma^1_1$-reflecting ordinal.
* It is unclear whether for non-Gandy $\alpha<\Omega$, $C(\Omega, \alpha)$ should canonically be $\alpha^+$ (the least admissible above $\alpha$), or the supremum of $\alpha$-recursive well-orderings (probably the latter), which (by definition) is $<\alpha^+$ for non-Gandy $\alpha$, and whether/how the additional structure below $\alpha^+$ affects ordinal assignments. Also, $\alpha^++1$-stable $\alpha$ are non-Gandy.

**Additional Properties**

An ordinal below $\Omega$ has degree a+$\Omega$ iff it is not above a certain ordinal and has degree *a*, or it is an admissible limit or a limit of admissible limits of ordinals of degree *a*. An ordinal above $\omega$ and below $\Omega$ in the notation is admissible iff its maximum degree has effective cofinality $\Omega$ using **O** (with effective cofinality of $\Omega$ being $\Omega$); equivalently iff it is a limit ordinal and every increasing sequence of ordinals in the notation system having that ordinal as a limit uses arbitrarily high ordinals up to the limit of **O** in the representation of the ordinals.

In the case **O** is based on $\Omega$ and C, a natural construction order for representable ordinals is the following: Start with 0, and then iteratively add the least ordinal that equals C(a, b) for *a* and *b* already constructed, with Ω added at stage $\varepsilon_0$. The following properties should hold. The construction order for representable ordinals below C(Ω, 0) agrees with their ordinal order. The construction order would remain unchanged if when adding C(a, b), we required that C(a, b) is in the standard form. The order-type of the construction order equals the order type of representable ordinals below C(Ω, 0). The construction order can also be relativized by starting with all ordinals (or all representable ordinals) below a particular one; such construction order should still satisfy analogues of the above properties.

Note that we are not defining C(a, b) when both *a* and C(a, b) are in a gap in the notation. One possibility is to allow ordinals not representable in the notation from below (that is using lesser ordinals as constants) to have the largest possible degrees. However, in that case some ordinals below C(Ω, 0) would have every degree below Ω (which contradicts the general definition of C). Another possibility is to use such *a* to fill in the gaps in the notation.

The notation system can just as well be defined above an arbitrary ordinal. Since **O** is a parameter, we can "stack" an ordinal notation on top of "itself" any finite number of times. However, stacking the notation system on top of itself an infinite number of times would lead to an ill-founded system (in a well-founded system, the order type using **O** and ordinals <a (given a<Ω and not counting ordinals ≥Ω) must be <Ω).

As noted previously, using $C_1$ rather than C can improve readability. Here is an example using $C_1$:

$C_1(\Omega^\Omega+\Omega^{\Omega^2}*3+\Omega^\Omega*4+\Omega^2*5+\Omega*6)$ is the 6th admissible ordinal after the 5th recursively inaccessible ordinal after the 4th recursively-Mahlo ordinal after the 3rd recursively hyper-Mahlo ordinal after the 1st $\Pi_3$ reflecting ordinal.

In $C_1(...+\Omega^a*b+...)$, $\Omega^a+b$ corresponds to going up by b (b is count) a-recursively inaccessible ordinals and their limits; if b is above the result, it encodes how much diagonalization to take (and similarly with other ordinals above the result). If a = $...+\Omega^c*d+...$, c encodes degree of recursive Mahloness and d encodes the number of times to take a limit operation where at successor steps the limits must correspond to c (where c encodes degree of recursive Mahloness). For example, a=$\Omega^2+\Omega*2$ corresponds to recursively Mahlo limit of recursively Mahlo limits of recursively hyper-Mahlo ordinals. Limit values of a,b,c,d are handled using the limiting process; for example $C_1(\Omega^{\Omega^\omega})$ is the least limit of recursively n-Mahlo ordinals for all finite n (and hence is not admissible).

A natural extension of the system with the same strength is to use $C_1$ instead of C in representing ordinals <Ω that are fix-points of x→$\omega^x$. For example, in testing for standard form, in $C_1(\Omega+d)$, d<Ω would be parsed for ordinals >d that are below Ω and not in the scope of an ordinal <$C_1(\Omega+d)$.

## 4.3 Assignment of Degrees

In this section, we present the conjectured canonical assignment of notations to ordinals. We assume that **O** uses CNF base $\Omega$ for **O**, or is or extends Veblen normal form (with exponentiation base $\Omega$). Note that the does not prevent an arbitrarily strong **O**.

In its most general form, given a sufficiently strong **O** for an appropriate theory T, the assignment can be used to represent canonical ordinals for T between an ordinal *a* (at least for *a* below non-Gandy ordinals or otherwise avoiding them) and the least $\Sigma^1_1$-reflecting ordinal >a. Such **O** would typically correspond to the proof ordinal of T. For weaker **O**, we get up to iterations of $\Pi^1_1$-reflection that are within **O**. The assignment works above an arbitrary ordinal *a*, using lower ordinals as constants, and allows choosing representations such that we need not worry about comparison of C-terms below *a*. For non-Gandy *a*, the definition works, but likely misses some canonical ordinals.

In formulating the assignment, our goal is that for appropriate T ranging from $\Pi^1_1$-$CA_0$ (or less) to (beyond) KP + $\forall \alpha \; \exists (\kappa^+$-stable $\kappa > \alpha))$, T + "every ordinal is assigned a term" is $\Pi^1_1$ conservative over T. For stronger theories using a stronger **O**, this should apply to every ordinal in the interval we cover. So our assignment has to work even in appropriate T too weak to prove well-foundedness, and we cannot require closure under subterms (consistently with T, a collapse of a fictitious ordinal might be well-founded).

**General Framework**

The assignment of ordinals assumes that the notation system has gaps; otherwise, all ordinals would be recursive; in the canonical assignment, gaps are set in the canonical way. Pick any sufficiently closed ordinal $\Omega$, say the least $\Sigma^1_1$-reflecting ordinal above the ordinals of interest (if building a notation system above an ordinal). For a simple suboptimal variation (or for the suboptimal old version), if using C or CNF base $\Omega$ for **O**, it suffices that $\Omega$ is $\varepsilon_{\Omega^++1}$-stable, and $\Omega^{++}$-stable should suffice for reasonable **O**. Alternatively, treat $\Omega$ as a syntactic construct and treat an ordinal $\geq \Omega$ either syntactically or as a pair ($\Omega$, ordinal using $\Omega$). We will define C(a,b) when *a* is representable using ordinals <C(a,b) as constants. This suffices for the system and allows extensions. In the system using C or CNF base $\Omega$ for **O**, the least non-representable ordinal is $C(C(\varepsilon_{\Omega+1},0),0)$; also, note that if $b \geq \Omega$, $C(a,b)=b+\omega^a$.

C(a,b) is the least ordinal above b that has degree *a*. Thus, it suffices to define degrees. A complication to recursively defining degrees is that notations can use higher ordinals for which degrees have not yet been constructed. We work around this by defining when:
x has degree d where d is a term that uses C, **O**, and ordinals <x as constants.

The recursion will assume that degrees for ordinals <x have been constructed. Therefore, (using standard comparison for C from "A Framework for Ordinal Notations" and the criterion for maximality of *a* in C(a,b)) we have comparison for terms that use C, **O**, and ordinals <x as constants (except possibly in weak theories). In turn, the degree of x lets us compute whether C(a,b)>x. Also, note that degree d makes sense for all ordinals y<x where y is above all constants used in d.

Every ordinal has degree 0. Assume that x is limit; otherwise x only has degree 0. Let us say that a degree (or (defined below) n-degree) d is <x maximal if for all sufficiently large y<x, d is maximal in C(d, y) (so C(d+1,y) > C(d,y)) (equivalently if there is y<x such that y is above all ordinals <x included in d and d is maximal in C(d, y)).

If we are promised that every C(a',b')≥x, then verification of <x maximality and comparison of such degrees is standard, depending only on comparison of constants <x. Also, <x maximality is defined/tested without referencing ordinals above x. For weak theories, we need an additional trick. Instead of using a single notation for d, in C(a',b') (b'<x) we can allow arbitrary b'<x that are at least as large as the standard choice (and referencing d only needs a single representative). We will use this to depend only on ordinals close enough to x. Moreover, if all b'<x (in C(..,b')) are identical, then the required testing for C(a',b')≥x depends only on the degree of x (with lower degrees permitting more forms, so underestimating a degree does not invalidate verification of having it, and the construction should work).

Note the similarity with [2.5 Reflection Configurations](#) since a degree gives a lower bound on the configuration of x above <x (if it exists).

x has degree d iff has x has <x-maximal degree d'≥d. (This holds because d does not use ordinals that are not representable from <x.) From now on, we mostly focus on <x-maximal degrees. In "x has <x-maximal d if ..", we assume that d is <x-maximal.

x has <x-maximal degree d+1 iff x is a limit of ordinals of degree d. This implies that x has degree d+1 iff x has degree d and (1) x is a limit of ordinals of degree d or (2) d is not <x maximal.

If d is a limit and the **O**-cofinality of d is <Ω, then x has degree d iff it has every degree <d.
*Notes:*
\* The cofinality is required to be <Ω because our notation system has gaps.
\* Weak theories typically will not prove that degrees <d are well-founded. However, the operative conditions for successor and **O**-cofinality Ω <x-maximal d' are phrased in terms of lower ordinals, which prevents infinite regress. An equivalent explicit definition for <x-maximal d (not of **O**-cofinality Ω) is that there are arbitrarily high d'<d such that x is a limit of ordinals of degree d'.

**Cofinality Ω (using CNF base Ω)**

If d has cofinality $\Omega$ in **O** (using CNF base $\Omega$ for **O**), separate out the rightmost 1 using CNF base $\Omega$ of d, and let n be its height (with the base being at height 1). Thus, n=2 for $..+\Omega^{..+1}$ (including $..+\Omega^1$), and n=3 for $..+\Omega^{..+\Omega^{..+1}}$, and so on.
x has degree d iff it has n-degree d or is a limit of such ordinals. n-degree is based on $\Pi_n$ reflection (and 1-degree is just degree), with the numbering of degrees chosen to be consistent between different n.
An ordinal x has n-degree $e+\Omega\hat{}\Omega\hat{}..\hat{}\Omega$ (using n-1 $\Omega$; using CNF; <x maximal degree) iff it is a $\Pi_n$ reflecting limit of ordinals of degree e. (Also, this is a special case of the next sentence since "$\Pi_1$-reflecting on" is the same as "limit of", except that only here do we allow e=0.)
More generally, x has n-degree $..(e+\Omega\hat{}\Omega\hat{}..\hat{}1)$ (using CNF; <x maximal degree) where the 1 is at height n and the '+' is at height m iff x is $\Pi_n$ reflecting and is $\Pi_m$ reflecting onto $..(e)$.
If $f = ..(e)$ has **O**-cofinality $\Omega$, $\Pi_m$ reflecting onto f means x is $\Pi_m$ reflecting onto ordinals of n'-degree f, where n' is the height of the rightmost 1 in f.
Otherwise, it means there are arbitrarily large <x-maximal g<f with **O**-cofinality $\Omega$ such that x is $\Pi_m$ reflecting onto g, with n' (the height of the rightmost 1 in g) is ≥m.

**Levels of Stability**

The above analysis also applies to a-stable x for a<x. We will use Veblen normal form base $\Omega$, for convenience using $\varphi_a(b)$ for Veblen $\varphi_a(\Omega+b)$ (a<$\Omega$), and $\varphi_0(b)$ for $\Omega^b$.

The height of the rightmost '1' in $d = ..+\varphi_{a_1}(..+\varphi_{a_2}(...+\varphi_{a_n}(..+1))..)$ (with $a_i<\Omega$) is $\omega^{a_1}+\omega^{a_2}+...+\omega^{a_n}$, but with 1 added if finite, and 1 subtracted (on the right) if the '1' is in a subterm of the form $\varphi_a(\Omega^1)$ or $\varphi_a(..+\Omega^1)$ (a>0 in both forms). (Note that ordinals like $\varepsilon_{\Omega+1}$ do not have cofinality $\Omega$.) The height of the rightmost '+' (above) is determined the same way, but without the subtraction.

Going further, the least κ-stable κ is $C(\varphi_{\varphi_\Omega(1)}(1))$ (with $\varphi$ as above). Diagonalization presents additional complications, but up to $\kappa^+$, we can address them as follows. The height (of '1' and '+' above) will still use $\omega^{a_1}+\omega^{a_2}+...+\omega^{a_n}$ (adjusted as above), but $a_i$ might be above x, in which case $a_i$ is converted into the order type of {h<$a_i$: h is <x-maximal}. (Also, in a weak theory, if this is ill-founded, we give up setting this degree, except through higher degrees.) Further still, $\Pi_\Omega$ reflection corresponds to κ being $\kappa^+$-stable (and to $\Pi^1_1$-reflection).

*Notes on transfinite reflection levels*:
\* For a>0, $\Pi_{\omega*a}$ reflecting means a-stable, and $\Pi_{\omega*a+n}$ reflection of x is defined using reflection of $L_{x+a}$ onto an ordinal below x for the appropriate set of formulas. The definition for reflection onto a class of ordinals is analogous.
\* Alternatively, transfinite $\Pi_a$ reflection (for a<x) can be defined using infinitary

formulas (with the set coding the formula in $L_x$). Note that every $\Pi_a$ formula for limit $a$ is $\Pi_{a'}$ for some a'<a.
* If *a* is an x-recursive well-ordering (potentially ≥x), $\Pi_a$ reflection can also be defined using almost self-referential formulas, with ψ(c,x) (c∈a), allowed to use ψ(c',y) with c'<c as a subformula (except that if c+ω>a, we have to count the quantifiers to confirm that the depth is ≤a).

## Iterations of $\Pi^1_1$-reflection

Starting with $\Pi^1_1$-reflection, the picture changes qualitatively because $\Pi^1_1$-reflection for κ can be iterated up to (informally) $κ^{++}$ times as opposed to up to $κ^+$ times for lower reflection levels, and iterations of the reflection capture the notation system. Essentially, ordinals $<κ^+$ can be used freely, because if a and b are two different recursive representations of α above κ (using constants <λ), then above λ<κ, that a and b denote the same ordinal is $\Pi^1_1$ and can be added as a reflection condition, so it does not matter whether the reflection is iterated along a or b. Now given an ordinal notation system **O** above a generic ordinal, interpret **O** above $κ^+$, and given an **O**-term A>0, Ath iteration of $\Pi^1_1$-reflection means being $\Pi^1_1$-reflecting onto Bth iteration for every B<A, with the 1st iteration simply means being $\Pi^1_1$-reflecting, or (as appropriate) $\Pi^1_1$-reflecting onto a given set.

Now, in the notation system (with x being $\Pi^1_1$-reflecting), parse the rightmost term of d until reaching either $φ_Ω(a)$ (a>0) or $a=φ_Ω(a)>Ω$. Collapse the terms (<Ω) in **O** representation of *a* into ordinals below $x^+$ using $a_i →$ order type of {h<$a_i$: h is <x-maximal}, and convert a into a' using **O** above $x^+$. Subtract 1 from a' (and store the result as a') if a' was b+n where n is finite and b is limit with **O**-cofinality <$x^+$. Now, having <x-maximal degree d means being a'-$\Pi^1_1$-reflecting and also being $\Pi_m$ reflecting onto ..(e) (or is a limit of such ordinals), where m and ..(e) are as above (see "Cofinality Ω (using CNF base Ω)" and its application to levels of stability; the rightmost '+' considered must be outside of a, so m<Ω).

Note that we do not require **O** to be well-founded, and d having **O** cofinality Ω is equivalent to every converging (to d) sequence having ordinals/terms unbounded in Ω. Assuming we have not missed anything, by plugging in appropriate **O**, we get the complete assignment below the least $Σ^1_1$-reflecting ordinal even in strong theories. For a sufficiently closed appropriate theory (such as $Z_2$), **O** would correspond to its proof-ordinal built above a generic ordinal, which might be available from other systems in the paper. Also, by using order-preserving bijections, the canonical assignment also applies to other appropriate ordinal notation systems.

*A simple suboptimal assignment:* There is also a relatively simple suboptimal assignment (that might not work in weak base theories) that should have the right limit when using Veblen normal form. As above, let e-$\Pi^1_1$ reflecting mean e-fold iteration of being $\Pi^1_1$-reflecting. C(0,x) will be the next admissible ordinal above x. x has <x maximal degree f iff x is f'-$\Pi^1_1$ reflecting (or even simpler, but with the wrong limit, f'-stable), or is a limit of such ordinals, where f' is roughly analogous to the old version below. Specifically, let $f_1$ be the order type of {h < f : h is <x maximal }. Let y be the order type of representable ordinals below C($\Omega$, x) (allowing x and ordinals <x as constants). Express $f_1$ using **O** base y and then change to **O** base $x^+$ (the least admissible ordinal above x) to get f'.

**Old Version**

*Notes:*
* The text branches from the main text starting with "Cofinality $\Omega$ (using CNF base $\Omega$)".
* The assignment below is retained for historical purposes but it misses some canonical ordinals starting with the least $\Pi_3$ reflecting ordinal that is $\Pi_2$ reflecting onto $\Pi_3$ reflecting ordinals, which leads it use higher ordinals (for given terms) than the canonical assignment. The canonical assignment only uses ordinals through $\Pi_n$ reflection. I would like to thank a reviewer with alias Hyp cos for pointing out the issue (in a related system) and suggesting the right strengths, which was instrumental in getting the new version above.

Otherwise, the fundamental sequence for d using CNF base $\Omega$ has length $\Omega$ and d is <x maximal. d has the form e + $\Omega^f$ (additive decomposition) where f is successor or is ≥$\Omega$. x has degree d iff x has admissibility degree f' and is a limit of ordinals of degree e, or a limit of such ordinals (that is a limit of ordinals of admissibility degree f' that are limits of ordinals of degree e), where f' is such that $f_{deg}$(f')=$\Omega^f$ (so f' is f or f-1; $f_{deg}$ is defined below). This reduces assignment of degrees (below $\varepsilon_{\Omega+1}$) to that of admissibility degrees. (Note: The notion of admissibility degrees used here corresponds to degrees of recursive Mahloness and beyond and not degrees of recursive inaccessibility.)

We now define admissibility degrees when the admissibility degree d is <$\varepsilon_{\Omega+1}$, which suffices for the notation system. Every ordinal has admissibility degree 0. Let $f_{deg}$ be the function converting admissibility degrees to degrees: $f_{deg}$(d) = $\Omega^d$, except that if d is a limit with fundamental sequence (using CNF base $\Omega$) shorter than $\Omega$, or if d is such a limit plus a finite ordinal, then $f_{deg}$(d)=$\Omega^{d+1}$. d has the form (CNF base $\Omega$) e + $\Omega^f$*g. x has admissibility degree d iff
Case A: $f_{deg}$(d) is not <x maximal. x has admissibility degree d' for the least d'>d such that $f_{deg}$(d') is <x maximal.
Case B: g is limit and $f_{deg}$(d) is <x maximal. x has admissibility degree e + $\Omega^f$*h for every representable h<g (equivalently, for cofinally many h<g).
Case C: g is successor and $f_{deg}$(d) is <x maximal. x is $\Pi_{2+f'}$ reflecting onto ordinals of admissibility degree e + $\Omega^f$*(g-1) (g-1 is the predecessor of g) where f' is derived from f as follows:

- Let $f_1$ be the order type of {h: $f_{deg}(e + \Omega^f*(g-1)+\Omega^h)$ is <x maximal and h is a representable ordinal <f} (so f' = $f_1$ = f if f≤x)
- Let y be the order type of representable ordinals below $C(\Omega, x)$. Express $f_1$ using CNF base y and then change CNF to base $x^+$ (the least admissible ordinal above x) to get f' (so f' = $f_1$ if f<$\Omega$).

If the notation system is restricted to terms below $\Omega^\Omega$, the assignment of gaps simplifies as follows: c = C(a,b) is the least ordinal consistent with its degree of recursive Mahloness and assignment of notations below c (here, the system is simultaneously built above all ordinals; the least notation incorrectly assigned would be $C(\Omega^{C(\Omega^\Omega, 0)+2}, 0)$). For every ordinal c=C(a,b) representable using lower ordinals as constants (even without the restriction on terms), c is admissible iff CNF base $\Omega$ fundamental sequence of *a* has length $\Omega$, and c is recursively f-Mahlo iff it is admissible and *a* is of the form ...+$\Omega^d$*e (CNF base $\Omega$) where d>f and f≤a; recursively 0-Mahlo means admissible.

## 4.4 A Step towards Second Order Arithmetic

The notation system below is only slightly more powerful, yet by introducing the idea of "correctness", it arguably brings us significantly closer to full second order arithmetic. Originally, the notation system was my attempt to reach the full second order arithmetic, yet the maximality condition used here gives a much weaker strength. Thus, while retained for historical purposs, the specific system here is essentially superceded by the main system. However, the "Levels of Correctness" part applies generally (including to the main notation system).

**Levels of Correctness**

In the above system (Degrees of Reflection), the reflection properties of C(a,b) (b<$\Omega$) change qualitatively as *a* crosses $\Omega$, and $\Omega$ itself (corresponding to $\Omega_2$ below) is qualitatively different. To reach second order arithmetic, we need $\omega$ levels like this (corresponding to $\Sigma_n$ correctness), which we abstract into the following general definition of correctness.

*Definition:* Let C be as in the General Notation. Let us say that a definition of correctness is compatible with C iff:
* Every ordinal has correctness 0; 0 only has correctness 0.
* The set of ordinals of correctness i is closed (i∈$\mathbb{N}$).
* If i<j, then correctness j implies correctness i.
* If i<j, the least ordinal above an ordinal b of correctness i is below the least such of correctness j. However, it is permitted that such ordinals are above all representable ordinals, and to define correctness only for representable ordinals.
* C(a,b) has correctness i>0 iff there is c with correctness i+1 with b<c≤a.
* If a is the least ordinal of correctness i above b, then C(a,b) < C(a+1,b).

For a notation system, let us use C, 0, $\Omega_i$ (0<i<$\omega$), with $\Omega_i$ being the least ordinal of correctness i. At this point (if we do not use gaps), despite its generality, C(a,b) is fixed apart from when *a* is maximal in C(a,b). Moreover, if using the strict version of C, the comparison of valid terms is independent of the notation system, and so is the determination of the correctness level of a valid term. For standard form, consistent with the main notation system below, we can use minimal b in C(a,b) and $\Omega_i$ in place of $C(\Omega_{i+1},0)$.

*Properties of the system*:
The least ordinal of correctness i>0 above b, $\Omega_i(b) = f^{(n-i)}(\Omega_n)$ whenever i≤n and $\Omega_n$>b and f(x)=C(x, b).
$0 < \Omega_1 < \Omega_2 < ...$
C(a, b) < $\Omega_i$ iff b < $\Omega_i$ and a < $\Omega_{i+1}$.
C(a, b) = $\Omega_i$ iff b < $\Omega_i$ and a = $\Omega_{i+1}$.
$\Omega_i$ is always maximal in $C(\Omega_i, b)$ (and as always, 0 is maximal in C(0, b)).
If an ordinal *a* has a positive but not infinite correctness n, b < a, and d is less than the least ordinal above *a* of correctness n+1, then C(d, b) < a. (Also, infinite correctness level can only occur in extensions of the system.)

*Reflection configurations:* (See [2.5 Reflection Configurations](#) for definition.) Let $\Omega_i(x)$ denote the least ordinal of correctness i above x. When defining a reflection configuration above b (or above <b), a reasonable notation is to replace $\Omega_i$>b (or $\Omega_i$≥b respectively) with $\Omega_i(x)$. This way, all reflection configurations are valid without an upper bound (unless we use a definition that breaks them), and we do not have to specify the domain explicitly. The comparison of configurations is lexicographical (in postfix form) after converting everything to the same $\Omega_i(0)$ and $\Omega_i(x)$ using $\Omega_j(x) = C(\Omega_{j+1}(x),x)$. In this notation, $\lambda x.\Omega_1(x) > \lambda x.\Omega_2(0)$, and the latter configuration is only valid for x>$\Omega_2$ (and 0 serves to avoid confusion).

**A Particular System**

Using the above template, let us define a system analogously to Degrees of Reflection, but with an additional built-from-below condition to ensure well-foundedness.

*Definition of the system:* If C(a, b) has maximum correctness n>0, let $\Omega$ be the least ordinal above b of correctness n+1, and $\Omega'$ the least ordinal above b of correctness n+2. Thus, we have $\Omega \le a < \Omega'$. For maximality of *a*, we require that the standard representation of *a* does not use ordinals that are above *a* but below $\Omega'$ except in the scope of an ordinal less than $\Omega$. In addition, as in the previous section, parse "*a*" from the root to branches until a constant or an ordinal below $\Omega$ is reached on every branch (at this stage, do not parse the ordinals below $\Omega$). For every such ordinal d < $\Omega$, we require that its standard form does not use ordinals strictly between d and $\Omega$ except in the scope of an ordinal less than C(a, b). If C(a, b) only has correctness 0, we simply require that the standard form of *a* does not use ordinals strictly between *a* and $\Omega$ except in the scope of an ordinal less than C(a, b), where $\Omega$ is the least ordinal above b of correctness 2.

This leads to a polynomial time algorithm for checking whether a particular representation is standard, and I expect there is also a polynomial time algorithm for converting arbitrary C-terms to the standard form.

For the conjectured canonical assignment:
Every ordinal has correctness 0.
An ordinal *a* has correctness 1 if it is admissible $>\omega$ (or a limit of such ordinals).
An ordinal *a* has correctness of n+1 (n>0) if it is a $\Sigma^1_1$-reflecting limit of ordinals of correctness n (or a limit of such ordinals).
*Note:*
A more natural assignment of correctness (corresponding to second order arithmetic) is given in the next section.

The notation system is conjectured to reach KP + {there is an ordinal of correctness n}$_n$. This is weaker than even KP + there is *a* that is $a^+$+1-stable, which is weaker than $\Pi^1_1$-CA$_0$ + lightface $\Pi^1_2$ comprehension. I have not ruled out that due to non-Gandy ordinals, the assignments should be lowered.

**Examples:**
If *a* has correctness n>0 (but not infinite correctness) and b is the least ordinal above *a* of correctness n and c < a, then $C(\varepsilon_{a+1}, c) = C(b, c)$.
$C(\Omega_{n+1}*2, 0)$ is the least admissible limit of ordinals of correctness n.
$C(\Omega_{n+1}^2, 0)$ is the least recursively Mahlo limit of ordinals of correctness n.
If *a* is the least ordinal above b of correctness n+1, then C(a+a, b) is the least admissible limit of ordinals of correctness n above b, and analogously with other large ordinal properties.

# 5 Main Ordinal Notation System

**Note:** The original title was "Ordinal Notation System for Second Order Arithmetic" and the system here is plausibly at least as strong as that.

## 5.1 Definition and Basic Properties

**Definition:** An ordinal *a* is 0-built from below from <b iff *a*<b
*a* is n+1-built from below from <b iff the standard representation of *a* does not use ordinals above *a* except in the scope of an ordinal n-built from below from <b.
(Note: "in the scope of" means "as a subterm of".)

The n$^{th}$ (n is a positive integer) ordinal notation system is defined as follows.
**Syntax:** Two constants (0, $\Omega_n$) and a binary function C.
**Comparison:** For ordinals in the standard representation written in the postfix form, the comparison is done in the lexicographical order where 'C' < '0' < '$\Omega_n$':
For example, C(C(0,0),0) < C($\Omega_n$,0) because 000CC < 0$\Omega_n$C.
**Standard Form:**
0, $\Omega_n$ are standard

"C(a, b)" is standard iff
1. "a" and "b" are standard,
2. b is 0 or $\Omega_n$ or "C(c, d)" with a≤c, and
3. *a* is n-built from below from <C(a,b) (use standard comparison to check).

**Notes:**
* For basic properties, see 2 A Framework for Ordinal Notations above. For example, C(a,b) = b+$\omega^a$ for *a* below the least $\varepsilon_x$>b.
* For b < $\Omega_n$, note that $\Omega_n$ is 1-built from below from b but not 0 built from below from b.
* A variation is to require *a* to be n-built from below from ≤b, which would simplify the definition a bit but lead to arbitrary restrictions like disallowing C(a,b) for a=$\Omega_2$+C($\Omega_2$*2+**C(C($\Omega_2$*3,0),0)**,0) and b=0 (the highlighted term is between b and C(a,b)).
* *Formalization of n-built from below:* As before, let $T_a$ be the parse tree of 'a': $T_a$ is the set of subterms of 'a', and for x and y in $T_a$, x⊏y means that x is a proper subterm of y; identical terms at different positions are distinguished (in the quantifiers and in ⊏). If a<$_n$b means *a* is n-built from below from <b, then a<$_0$b ⇔ a<b and a<$_{n+1}$b ⇔ ∀x∈$T_a$ (x>a) ∃y⊒x y<$_n$b. For example, a<$_2$b ⇔ ∀x∈$T_a$ (x>a) ∃y⊒x ∀z⊏y (z>y) ∃t⊒z t<b.

For n=0, this system just reaches $\varepsilon_1$, for n=1, this is the notation system for 2.3 Bachmann-Howard ordinal given earlier, and for n=2, the notation system is superficially similar to (but perhaps much stronger than) the one in 4 Degrees of Reflection. A natural conjecture (see 5.2 Strength of the Main System) that the strength of the n[th] ordinal notation system is between $\Pi^1_{n-1}$-CA and $\Pi^1_n$-CA$_0$ (or higher), and thus the sum of the order types of these ordinal notation systems is the proof-theoretical ordinal of second order arithmetic (or higher).

**Combined System:** The ordinal notation systems are best combined into one system as follows:
Constants 0 and $\Omega_i$ (for every positive integer i), and binary function C.
$\Omega_i$ = C($\Omega_{i+1}$, 0) and the standard form always uses $\Omega_i$ instead of C($\Omega_{i+1}$, 0).
To check for standard form and compare ordinals use $\Omega_i$ = C($\Omega_{i+1}$, 0) to convert each $\Omega$ to $\Omega_n$ for a single positive integer n (it does not matter which n) and then use the n[th] ordinal notation system.

**Nonstandard forms:** To make C a total function for *a* and *b* in the notation system, let C(a, b) be the least ordinal of degree ≥a above b, where the degree of $\Omega_i$ is $\Omega_{i+1}$ and the degree of C(c,d) is c if "C(c,d)" is the standard form. I believe that this is recursive with the conversion to standard form as follows.

**Conjectured conversion to the standard form for the n[th] notation system:**
To convert 'C(a, b)' to the standard form, first convert 'a' and 'b'. Next, recursively minimize b by replacing it with d for as long as b is C(c, d) and c<a. If *a* is not n-built from below from b, then perform an in-order right-to-left traversal of the term tree of *a* to find the first counterexample: a < $a_1$ < $a_2$ ... < $a_n$ ($a_{i+1}$ is a subterm of $a_i$). This occurrence of $a_n$ is part of $C(a_n, d)$ with $C(a_n, d) < \Omega_n$. Replace this $C(a_n, d)$ with $\Omega_n$, and delete everything to the left of this $\Omega_n$, and add the right number of 'C(' to make the new 'a' a valid term, and convert the new 'a' to the standard form.
**Note:** By converting between different $\Omega_i$, the conversion for the full notation system follows.

**Proposition:** In the *n*th system, given a standard term *a* (in prefix form), every nonempty suffix of *a* augmented with the corresponding number of 'C' on the left is standard.
*Note:* For example, for standard C(C(a,b),c), C(b,c) is standard. Also, in the combined system, the result is standard after a possibly repeated conversion $C(\Omega_{i+1},0) \to \Omega_i$.
*Proof:* By recursion, it suffices to check when a=C(c,d), c is modified into c', and c' is standard. Since c' has no subterms between c' and c (or between C(c',d) and C(c,d), and similarly for certain subterms of c') and since the subterm tree is sufficiently preserved, the n-built from below criterion will carry over to c'.

**Variations using one-variable C:** If the system is expressed using $C_1$ (one variable C), the maximality condition on *a* in $C_1(a)$ may not be the simplest. A natural variation (for the $n^{th}$ notation system) is to require *a* to be n-built from below from <$C_1(a)$ (terms use 0, $\Omega_n$,'+', x→$\omega^x$, and $C_1$ for fix-points of x→$\omega^x$) and another natural variation (though it is unclear whether it is too restrictive to get the strength) is to simply require for *a* to be n-built from below.

**Passthrough extension:** It is natural to try to extend the system by combining it with ideas from [8 Built-from-below with Passthrough](#). For example, let *a* be n+1-built from below from <b if the representation of *a* does not use ordinals above *a* except in the scope of an ordinal n-built from below from <b, where in representing C(c,d) for c<a, one only counts C(c,d) and ordinals ≤d and their subterms. This allows additional terms such as $C(\Omega_2+C(\Omega_2*2+C(\mathbf{C(\Omega_2*3,0)},C(\Omega_2*2,0)),0),0)$ (the term is C(a,0) and c is highlighted). In defining the extension, we have to be careful to prevent passthrough from wrapping terms into built-from-below terms, which could then be collapsed into ill-foundedness. To achieve this, in C(a,b) we require that *a* is a,n-built-from-below from <C(a,b), denoted a<$_{a,n}$C(a,b).
a<$_{a',0}$b ⇔ a<b
a<$_{a',n+1}$b ⇔ ∀x∈$T_a$ (x>a) ∃y⊒x (y<$_{a',n}$b ∨ y⊐x ∧ d(y)<a' ∧ ∀z (x⊑z⊑y) z≥y)
where d(C(e,f))=e (for standard C(e,f)) and d(y)=∞ otherwise. The combined system can be defined as above. A natural restriction is to prevent entering n-built-from-below stage while we are in the passthrough. The passthrough might make the system much stronger (especially if delayed diagonalization is not a problem here), but given our lack of knowledge, also has a high risk of ill-foundedness. If full passthrough is a problem, a restricted version can still likely be used for iterating n-built-from-below, as described in [7 Iteration of n-built from below](#).

**Reflection configurations:** There is also a possible variation on the main system using [2.5 Reflection Configurations](#) (see there for terminology); its strength and well-foundedness are unclear. The nth system will require that in C(a,b), r(a) is r(a),n-built-from-below, denoted by $B_{r(a),n}(r(a))$, or we are not using passthrough, 0,n-built-from-below. The combined system is as above.

$B_{a,0}(d) \Leftrightarrow d < \lambda x.x$ (i.e. reflection configuration d is a constant (below x))

$B_{a,n+1}(d) \Leftrightarrow \forall x \in T_d\ (x \sqsubset d \wedge r(x) > d)\ \exists y \sqsubset d\ (r(x) \subseteq r(y) \wedge (r(y) < a \vee B_{a,n}(r(y))))$.

*Notes:*

\* An equivalent construction that avoids nested configurations is to require $a \prec_{r(a),n} C(a,b)$ (using standard comparison for C(a,b)), or $\prec_{0,n}$ if not using passthrough, with

$d \prec_{a,0} b \Leftrightarrow d < b$ (note that d is an ordinal)

$d \prec_{a,n+1} b \Leftrightarrow \forall x \in T_d\ (x \sqsubset d \wedge r(x) > r(d,b))\ \exists y \sqsubset d\ (r(x) \subseteq r(y) \wedge (r(y) < a \vee y \prec_{a,n} \text{parent}(y)))$.

\* The scope of reflection configurations (which is narrower than the scope for ordinals) already gives us a form of passthrough, but the passthrough condition becomes relevant after our parsing leads into a higher configuration. For example, in the final stage, rather than being built-from-below (and thus in a sense, recursive), we have built-from-below with passthrough, which might be much broader.

\* Given the triple permissiveness (n-built-from-below, reflection configurations, passthrough), this is the place for the interested reader to look for ill-foundedness, but quite possibly this is also the only version that actually reaches second order arithmetic with projective determinacy.

## 5.2   Strength of the Main System

**Possibilities for the Strength**

Behind its simple definition, the main notation system hides subtle and unknown structure. Similar patterns repeat at different levels of the strength hierarchy, and here are five of the possibilities for its strength: (1) not well-founded, (2) rudimentary set theory + {there is $\alpha^{(+)^n}$-stable $\alpha$: $n \in \mathbb{N}$} (3) second order arithmetic, (4) rudimentary set theory + {there is n-subtle cardinal: $n \in \mathbb{N}$}, (5) second order arithmetic with projective determinacy (perhaps using canonical projective ordinals).

Without a proof or a clear understanding of its structure, we cannot rule out (1). At the moment, (2) represents a lower bound on our intuitive ability to embed descriptions of ordinals into the notation system (and even there, I did not properly analyze structure at non-Gandy ordinals, which start appearing before $\alpha^+ + 1$-stable); we do not expect the strength to be that low, but we do not yet see beyond that with a sufficient clarity. (3) corresponds to the similarity between the nth system and around n alternating quantifiers, and to a rough similarity between $C_i$ in the passthrough system ([8 Built-from-below with Passthrough](#) below) and the usage of $\Omega_i$. (4) often shows up in certain loosely similar problems (see Harvey Friedman's work and also (Taranovsky 2012)). Moreover, the relationship between the passthrough system and the main system at level of $\Pi_n$-reflection suggests a strength beyond (3). Going further, despite inclusion of small large cardinals, L has a $\Delta^1_2$ well-ordering (of ordinals countable in L). If second order quantifiers are fully used (corresponding to a $\Delta^1_n$ well-ordering), the strength would correspond to (5). A number of different natural additions to second arithmetic, including (conjecturally) comparability of Wadge degrees lead to projective determinacy. Our guess for the main system is (3), but we expect that a reasonable modification can be used to get (5). Perhaps, some combination of passthrough with n-built-from-below works (an example is given above). See [6.2 New Analysis](#) for additional details (and possibilities).

Also, an intuitive analysis by user Hyp cos of a fragment of the main system (unmodified version) with a long list of examples can be found at http://googology.wikia.com/wiki/User_blog:Hyp_cos/TON,_stable_ordinals,_and_my_array_notation ; I have not verified its correctness.

**Differences from Degrees of Reflection C**

The n=2 notation system in [5 Main Ordinal Notation System](#) is not identical to the one in [4 Degrees of Reflection](#). The notation system for degrees of recursive inaccessibility corresponds to the one in "Degrees of Reflection", and the assignments of ordinals in "Degrees of Recursive Inaccessiblity" and "Degrees of Reflection" should not be used for the main system. Let the standard C denote denote the notation system in "Main Ordinal Notation System".

$C(C(\Omega+C(\Omega*2,0),C(\Omega*2,0)),0)$ appears to be the "least" term valid for "Degrees of reflection" C but not for the standard C. However, this does not appear to adversely affect the strength of the standard C since enough ordinals remain available for the collapse. For example, in n=2 system, for every *a* built from below, for every b that uses only ordinals <a, C(C(a+b,0),0) is standard (after converting a+b to the standard form).

$C(\Omega_2+\varepsilon_{C(\Omega_2*2,0)+1},0)$ appears to be the "least" term valid for the standard C but not "Degrees of Reflection" C, and the difference appears to make the standard C for n=2 much stronger than C in "Degrees of reflection".

**Conjectured Hierarchy**

For standard C (with $\Omega=\Omega_2$), if $d=C(\Omega,C(\Omega*2,0))$, then $C(\Omega+d,0)$ should be the least recursively inaccessible ordinal. (Intuitively, since d is admissible but not limit of admissibles and for cofinally many representable e<d, C(1,e,0) is standard, C(1,d,0) should also be admissible.)
Continuing upward, $C(\Omega+d*2,0)$ should be the least recursively hyperinaccessible ordinal; $C(\Omega+d^2,0)$ -- least recursively Mahlo, $C(\Omega+d^d,0)$ -- least $\Pi_3$ reflecting, and so on, analogously to the description from "Degrees of Reflection".

The hierarchy would suggest that $C(\Omega*2, 0)$ should be the least ordinal x such that $L_x$ is a $\Sigma_1$ elementary substructure of L, and this hierarchy is formalized below in [8 Built-from-below with Passthrough](#), which is a variation on a fragment of the main system. However, while the passthrough system makes certain collapses free, 2-built-from-below (without passthrough) gives one free collapse, which might be insufficient as we approach an ordinal *a* that is $a^{++}$-stable (or just $\Sigma^1{}_1$-reflecting). Roughly speaking, getting an ordinal d corresponding to a large level $<a^{++}$ involves one collapse and arguments corresponding to $<a^+$, with a second collapse to get those arguments.

The level at which the main system reaches $a^{++}$ stable ordinals is unclear. Essentially, the worst case scenario is that correctness n>1 corresponds to the closure of $\{\alpha: \alpha \text{ is } \alpha^{(+)^n}\text{-stable}\}$ (or even less, with $\Omega_2$ possibly being the least $\Sigma^1{}_1$ reflecting ordinal). Consider $c = C(\Omega_2+a,b)$, and let f be an appropriate recursive ordinal function (using parameters $<c$; $f(x)>x$), and consider how f can be placed in *a*. In constructing c, the final collapse in the definition of f is free, which allows using f in many but not all places inside *a*. If we can show that the permitted places are enough, then $\alpha = C(\Omega_2*2,0)$ would already correspond to $\alpha^{++}$-iterations (or more if there is some stronger other structure). Note however, that places inside *a* that will be collapsed further impose a limit on the ordinals used in f (if placed there), which would not work for getting the strength.

While the definition of the main system is perhaps simpler than the passthrough system, the relationship between the ordinals and the notations appears more direct in the passthrough system, especially if the main system does not reach beyond second arithmetic. Assuming it is strong enough, the passthrough system probably works better for proving the exact ordinal strength of $Z_2$.

**Higher n:** If the strength of the main system corresponds to second order arithmetic, some conjectured assignments are as follows.
$\Omega_n$ (n>1) - as in the next section
$C(C(C(\Omega_3+1,0),0),0)$ - the proof ordinal of $\Pi^1{}_2\text{-CA}_0$, and analogously for $\Pi^1{}_n\text{-CA}_0$.
$C(C(C(\Omega_3+\varepsilon_0,0),0),0)$ - the proof ordinal of $\Delta^1{}_3\text{-CA}$, and analogously for $\Delta^1{}_n\text{-CA}$.
$C(C(C(\Omega_3*2_0,0),0),0)$ - the proof ordinal of $\Pi^1{}_2+\text{TR}_0$, and analogously for $\Pi^1{}_n+\text{TR}_0$.

**Levels of Correctness**

For terms in the notation system, assign correctness levels as in [4.4 A Step towards Second Order Arithmetic](#) above. (See above for basic properties, but for example, $\Omega_i$ (i>0) is the least ordinal of correctness i.) If the notation system corresponds to $Z_2$, then we expect the following under the conjectured canonical assignment:
\* An ordinal κ has correctness 1 iff it is admissible $>\omega$ or is a limit of admissible ordinals.

\* An ordinal κ has correctness n+1 with n>0 iff $L_κ$ is a $Σ_n$ elementary substructure of $L_ρ$.
**Note:** Here, ρ is an ordinal such that $L_ρ$ satisfies ZFC minus power set. For definiteness, let ρ be $ω_1^L$ since second order arithmetic does not prove existence ρ with $L_ρ$ satisfying ZFC\P, but it proves as a schema that $L_{ω_1^L}$ satisfies ZFC\P.

## 5.3 Old Analysis and Examples

**Note:** Below is a set of (mostly wrong) guesses (retained for historical purposes) that were developed before a more modern understanding of the properties of the main system. The assignments assume that the main system corresponds to second order arithmetic.

**Ordinal Assignments**
After $Π_n$ reflecting ordinals, the next set theoretical concept is that of stability. Here are ordinals corresponding to different levels of stability.
$C(Ω_2^+,0)$ -- the least ordinal *a* that is stable up to $ε_{a^++1}$ (note that $Ω_2^+$ equals $C(C(Ω_3,Ω_2),Ω_2)$)
$C(C(Ω_3,Ω_2),0)$ -- the least ordinal stable up to 2 admissible ordinals
$C(C(Ω_3,Ω_2)*Ω_2,0)$ -- the least ordinal that is $Π_2$ reflecting onto ordinals that are stable up to two admissible ordinals, and similarly with other levels of reflection
$C(C(Ω_3,C(Ω_3,Ω_2)),0)$ -- the least ordinal stable up to 3 admissible ordinals, and so on
$C(C(Ω_3+Ω_2,0),0)$ -- the least ordinal stable up to a recursively inaccessible ordinal
$C(C(Ω_3+Ω_2^2,0),0)$ -- the least ordinal stable up to a recursively Mahlo ordinal, and so on
$C(C(Ω_3+C(Ω_3,Ω_2),0),0)$ -- the least ordinal that is stable up to a larger ordinal that is stable up to two admissible ordinals
$C(C(C(Ω_3*2,0),0),0)$ -- proof theoretical ordinal of $Π^1_2\text{-CA}_0$
$C(C(Ω_3*2,0)+1,0)$ -- the least ordinal *a* such that $L_a$ models KP + $Σ_1$ separation (which has the same strength as $Π^1_2$-CA + TI), so its proof theoretical ordinal appears to be $C(C(C(Ω_3*2,0)+1,0)^+, 0)$.
$C(C(C(Ω_3*ε_0,0),0),0)$ may be the proof theoretical ordinal of $Δ^1_3$-CA, which has the same strength as KP + {there is a $Σ_1$ elementary chain of length *a*: $a<ε_0$}.

**Possible Alternative**
One issue with the notation system is that while $C(Ω_{n+1}+1, 0)$ (n>0) is the height of the least $Σ^1_n$-correct model of $Π^1_n$-CA, the collapse $C(C(Ω_3+1, 0), 0)$ corresponds to a weaker theory than $Π^1_2\text{-CA}_0$. The cause of this is our maximality condition for degrees which causes for example $C(ε_{Ω_2+1}, 0)$ to equal $C(Ω_2^+, 0)$. If we could devise a corresponding notation system that distinguishes such terms, here are guesses for some of the proof theoretical ordinals:

For n>0, $C(...C(\Omega_{n+1}+1, 0)..., 0)$ (with n+1 Cs) -- the ordinal for $\Pi^1_n\text{-CA}_0$.

$C(C(C(\Omega_3, C(\Omega_3+1, 0)), 0), 0)$ -- the ordinal for $\Pi^1_2\text{-CA} + \text{TI}$.

$C(\Omega_3*2, 0)$ -- the least ordinal κ of correctness 2 with $L_\kappa$ satisfying $KP + \Delta_2$ separation.

$C(C(C(\Omega_3*2, 0), 0), 0)$ -- the ordinal for $\Pi^1_2\text{-TR}_0$.

$C(C(C(\Omega_3, C(\Omega_3*2, 0)), 0), 0)$ -- the ordinal for $\Delta^1_3\text{-CA} + \text{TI}$.

# 6 Beyond Second Order Arithmetic

## 6.1 Old Analysis

**Notes:**
* This subsection is retained for both intuitive and historical value. A more precise analysis of how C might reach ZFC and beyond is in the next section. Also, the three variable C described here corresponds with [8 Built-from-below with Passthrough](#) system.
* Both the old (and especially) the new analysis assume some understanding of set theory. Useful references include (Jech 2006) and (Kanamori 2008).

To go beyond second order arithmetic, we need transfinitely many degrees of correctness. Cardinals will be ordinals that cannot be reached from below no matter how large the degree of correctness is. Let Ω be the least uncountable cardinal and b a countable ordinal, as computed in the model (specifically, in L). If $a < b$, then b having correctness $\omega*(a+1)$ may be defined as $L_{b+a}$ being elementarily embeddable $L_{\Omega+a}$. Correctness $\Omega+\omega$ may correspond to $L_{b+b}$ being elementarily embeddable in $L_{\Omega+\Omega}$, and similarly for $\Omega^2+\omega$ and $L_{\Omega*\Omega}$, and so on.

For a notation system, we can try to use a total ternary function C such that C(a, b, c) is the least ordinal above c of correctness *a* and for that correctness of degree b. If we treat (a, b) like an ordinal, then the function satisfies formal requirements of C (as described in the general notation), so the only issue is specifying when (a, b) is maximal. For *a*=0, the maximality of b is arguably like in [3 Degrees of Recursive Inaccessibility](#) (that is the standard form of b uses only ordinals below b), but the general case is unclear.

I do not have the conditions for *a* and *b* in the presence of uncountable cardinals. However, roughly speaking, the condition for *a* is being definable from below from <C(a,b,c). If a'<a, then C(a',b,c) is definable in $L_{\Omega+a}$ for every representable value of b provided that $c < \Omega$ and a' and c are definable in $L_{\Omega+a}$ (Ω+a (in both instances) is not sharp here, but it communicates the idea).

I also note that, for every finite n, to get an ordinal notation system for ZFC\P + "$\omega_n$ exists", it is sufficient to describe an ordinal notation system for ordinals between $\omega_n$ and $\omega_{n+1}$ where ordinals $\leq\omega_n$ are given as constants and ordinals $\geq\omega_{n+1}$ are represented in terms of ordinals $<\omega_{n+1}$ using a given notation system

**O**. One can then stack the resulting systems on top of each other any finite (but fixed) number of times to get the full system for a particular n.

If the definition of C(a,b,c) is worked out, one can extend the system by adding a function that enumerates L-cardinals (or just cardinals). To get inaccessible cardinals, one can use constant 0 and a 4-variable C(a,b,c,d) where *a* indicates degree of inaccessibility. One guess is that going beyond 4-variables corresponds to higher order set theory, and that the strength for finite variable C corresponds to ZFC + {n-ineffable: n<ω} (plus a conservative higher order set theory extension to extend Ord to match the notation system). One can go further by adding ordinal Ω and using higher ordinals in place of n-tuples. To go still further, one would use $Ω_1$, $Ω_2$, ..., (analogously to the notation system for second order arithmetic) with the hope/ideal to capture second order arithmetic plus projective determinacy. However, we are still far from that.

**Collapsing functions**

One approach to find ordinal notation systems for ZFC and beyond is to find a set of functions that is rich enough to capture the set-existence principles but tame enough for comparison to be recursive. Let us start with 0, ordinal addition, ordinal exponentiation base ω, and add aleph function that enumerates infinite cardinals as computed in L. A key property of uncountable cardinals is unreachability from below; for example, every infinite model has a countable submodel. So to capture the unreachability, let us add function f: f(a, b) → the least ordinal c such that c is not definable in $L_a$ using parameters ≤b.

Two questions to ask are:
Is the resulting ordinal representation system recursive? (That is, is there an algorithm for comparing terms?)
If yes, how strong is the resulting system? Does it capture the full strength of the underlying set theory or does it correspond to recursive analogue (of cardinals) or something in between?

If f is not sufficient, one can try a full collapsing mapping and instead of just its critical point:
Coll(a, b) → the collapsing mapping for least elementary submodel of $L_a$ that contains all ordinals ≤b (and in that case, include terms like Coll(a,b)(c) in the notation system)

To try to reach ZFC and beyond, we can add a 3-variable function
f: f(a, b, c) → the least a-inaccessible (in L) L-cardinal k > c that that is not definable in $L_b$ using parameters <k
and ask whether the resulting representation system is still recursive, and if so, what is its strength.

To (try to) go to levels of indescribability, one can add function
g: g(a, b, c) → the least cardinal d > c such that for some cardinal e>d, there is elementary embedding j:Coll(e, d)($L_e$) → $L_a$ with crit(j) = d and j(d) = b. To (try to) get to n-ineffable cardinals for all n<ω, one can use a predicate for n-reflective cardinals as computed in L (reflective cardinals are defined in (Taranovsky 2012)), and define $g_1$, $g_2$, ... that are like g except that Coll and j must also be elementary with respect to n-reflective cardinals (for n<i for $g_i$), but allowing Coll and j to "move" the predicates for reflective cardinals by application.
To go beyond that, one would use mice and the associated embeddings.

**An example notation system**

[*Update:* Retained for historical purposes, but 7 Iterations of n-built from below is (at present) a better candidate. Also, it is unclear whether the system is related to 8 Built-from-below with Passthrough, and the caveat of it being a guess (that possibly was not thought through) applies.]
To inspire future work, here is how a notation system (for ZFC+{n-ineffable}$_{n<\omega}$) might look like, though it would be lucky if this particular system is well-founded and has the right strength. The system uses constant 0 and multivariable C. $C(a_1, a_2, ..., a_n, b)$ corresponds to 2-variable $C((a_1, ..., a_n), b)$ if $(a_1, ..., a_n)$ is treated as an ordinal with lexicographical comparison, where after all leading 0s are removed, a longer sequence is larger than a shorter one, and $(a_1)=a_1$, and standard form prohibits $a_1=0$. Thus (as written in 2.2 Basic Properties section for general C), it suffices to state the maximality condition.
$(a_1, ..., a_n)$ is maximal in $C(a_1, ..., a_n, b)$ iff for each i 1≤i≤n, $a_i$ is <$(a_1, ..., a_n)$ built from below from <$C(a_1, ..., a_n, b)$.
*a* is <b built from below from <c if for every subterm $e_i$ of *a* where $e_i$ is is part of $C(e_1, ..., e_m, f)$ (as an immediate subterm in the term tree of a) and $e_i$>a and on the path from *a* to $e_i$ all ordinals are ≥c and ≤a, we have $(e_1,...,e_m) < b$ and $e_i$ is <b built from below from <c.
The hierarchy is a bit different from the notation systems in previous sections. If d=C(1,0,C(2,0,0)), the least recursively inaccessible ordinal should be C(1,d,0), least admissible limit of recursively inaccessibles should be C(1,d*2,0), least recursively Mahlo C(1,$d^2$,0), and so on.
Instead of n-tuples (as in $(a_1, ..., a_n)$), we can use a stronger notation system — the construction can be formally generalized into a mapping that turns ordinal notation systems above generic ordinals (notation system above a generic ordinal is defined earlier in the paper) into stronger ones. As a particular example, we can use 2-variable C (in place of multivariable C), 0, and a large ordinal Ω, and treat ordinals ≥Ω as syntactic constructs. (Formally, given *a*, let variables range over subterms of *a*, including their position in *a*, and let '⊑' denote subterm (not necessarily a proper subterm, but different positions in *a* are distinguished), and '<' compare variables as ordinals. The maximality condition for *a* in C(a,b) is ∀s,t (s<t<Ω ∧ t⊑s ∧ ¬∃u⊒s s<u<Ω ∧ ¬∃u⊒t t<u<Ω ∧ ¬∃u⊒t u<C(a,b) ⇒ ∀v⊒t (¬∃w<Ω (t⊑w⊑v ∧ w≠t ∧ w≠v) ⇒ v<a)). Also, b<Ω⇒C(a,b)<Ω.)

## 6.2   New Analysis

*Introduction and notes:*

* This section provides insights about ordinals and cardinals (especially in L) and how a C-like system is likely to work there. The specific analysis is for a system behaving like [8.2 Degrees of Reflection with Passthrough](#) but without delays in diagonalization. We refer to it as the n=2 system, though whether the main system for n=2 behaves like this is not clear. The system is (actually or hypothetically) extended to higher n analogously to the main notation system.
* If either the main system, or Degrees of Reflection with Passthrough, or a passthrough extension of the main system (or a system using reflection configurations) works like that, the analysis is (approximately) accurate for that system.
* The ordinal assignments may seem high, but keep in mind that in L, all reals are $\Delta^1_2$ in a countable ordinal, and if the notation system has mastered $\Sigma_\alpha$-elementary substructures in L (and transitive collapses), the constructible versions of (small) large cardinals start to resemble the recursive versions (but with additional structure).
* If Degrees of Reflection with Passthrough (or a fragment considered, and especially if not using reflection configurations) is too weak, we conjecture that the assignment can still be defined and formalized by relying on the extensibility of the system (while keeping the assignment unchanged), though it is unclear how the extensibility should be formalized. Such an assignment would still qualitatively enhance our understanding of ordinals and cardinals, but there are points where it would fade out. One possibility is that for a theory T, models of T at definability level α are only captured if T is not too strong relative to (and has an appropriate definition at) the next definability level above α. Thus, the least rank-initial model of ZFC (or $V_\delta^L \vDash$ ZFC) would likely be captured, but because ZFC is too strong relative to the concept of transitive models, the least transitive model of ZFC would not. Or more precisely (if the system is parameterized by an ordinal notation system **O** built above all ordinary ordinals), the supremum of the heights of the least transitive models of finite fragments of ZFC would then likely equal C(1,δ,0) where δ corresponds to the **O**-term (or diagonalization level) corresponding to the proof ordinal of ZFC.

**Note:** We will work in L (where V=L).

*Previous introduction:* To inspire future research, here is an optimistic guess about the strength of C. The assignments might be too low if I missed structure in the notation system or too high if I missed structure in L. I assume that the n=2 system is analogous to [8 Built-from-below with Passthrough](#) (and in particular Degrees of Reflection with Passthrough) and in effect analyze that system, and that section (below) will help the reader to understand the discussion here. However, whether the analogy holds is unclear, and should the main system be too weak, the comments below for the n=2 system might still work for Degrees of Reflection with Passthrough (though it could be too weak as well), with n=3 and higher systems not yet defined (but see passthrough extensions to the main system for a possibility).

C(Ω*a+b,c) for b<C(Ω*a+b,c) should be ordinal number $\omega^b$ of definability/reflection level *a* above c (here $\Omega$ is $\Omega_2$). Otherwise (b>C($\Omega$*a+b,c)), b has a canonical definition "b" (using ordinals <C($\Omega$*a+b,c) as constants), and C($\Omega$*a+b,c) is roughly the least ordinal corresponding to definability level *a* model of "b" (above c). For example, a=1 corresponds to transitive models, and (working in L and assuming we have enough strength and there is no acceleration or delay) C($\Omega+\omega_\omega^L$,0) should be the least ordinal d such that for every finite n, there is a transitive model of ZFC\P+"$\omega_n$ exists" below d.

For simplicity (and insight), let us analyze the $\Sigma_2$ hull of L. Here, we can have large cardinals but every set will be $\Sigma_2$ definable, or $\Sigma_1$ with the cardinality predicate. When *a* is sufficiently definable, c having definability/reflection level *a* corresponds to existence of structures <c of arbitrary complexity for (sufficiently definable) definability levels a'<a, allowing ordinals <c as parameters. For example, the least ordinal $\kappa$ with $L_\kappa <_{\Sigma_1} L$ is the least ordinal such that transitive models for arbitrary sufficiently satisfiable axioms exist below $\kappa$. Assuming the base theory gets its strength from the large cardinals it asserts, these axioms correspond to large ordinals/cardinals: structure for a transitive model → "b exists" → b (where large ordinal/cardinal *b* is canonically defined from below), and the models roughly correspond with C($\Omega$+b,0), so C($\Omega$*2,0)=$\kappa$.

Approximate definability level correspondence for c=C($\Omega$*a+$a_1$,b):
small 1+a -- $L_c$ is a $\Sigma_a$ elementary substructure of $L_{\omega_1}$ (c is countable, *a* may be transfinite, a=0 corresponds with $L_c$ being admissible or limit of admissibles).
large a<$\Omega$ -- c is not $\Sigma_1$ definable in ($L_{\kappa+a}$, ∈, Card) (Card is predicate for cardinals (in L)) allowing ordinals <c as constants, where $\kappa$ is the least cardinal above c.
a=$\Omega$ -- cardinals. C($\Omega^2$+a,0) is $\omega_{\omega^a}$ for *a* below the least fix-point of the aleph function.

To go further, let us list some cardinals (in L) in the increasing order below the least inaccessible cardinal $\lambda$ in L:
the least $\kappa$ such that $V_\kappa$ is a model of ZFC.
the least cardinal $\kappa$ such that for some $\alpha$, $L_\alpha$ satisfies "P holds and $\kappa$ is inaccessible" where P is some large cardinal axiom that is compatible with L.
the least $\kappa$ such that $V_\kappa$ and $V_\lambda$ satisfy the same $\Sigma_2$ statements.
the least $\kappa$ such that $V_\kappa <_{\Sigma_2} V_\lambda$. $\kappa$ is also the least fix-point of f where f(a) is the least b such that $V_b$ and $V_\lambda$ satisfy the same $\Sigma_2$ statements with parameters in $V_a$.
the least $\kappa$ such that $V_\kappa <_{\Sigma_a} V_\lambda$. For large 'a', one analog of this is the supremum of ordinals <$\lambda$ that are $\Sigma_1$(Card) definable in $L_a$ using parameters <$\kappa$.

Note the parallel between these notions and ordinals ≥$\lambda$ that are below $\lambda^+$ (the cardinal successor of $\lambda$), and our intuitive analysis suggests this analogy also applies to the notation system (hence the choice of the notation for the least inaccessible cardinal).

Let C(a,b,c,d) = C($\Omega^2$*a+$\Omega$*b+c,d) and C(b,c,d) = C(0,b,c,d), and $x^+$ = C(1,0,0,x)

(the cardinal successor of x), and let S=C(1,1,0,0).
C(1,0,S,0) -- the least fix-point of x→ω$_x$.
C(1,0,C(1,0,S),0) -- a lower bound on the least power-admissible ordinal; we have not analyzed whether it is even higher.
C(1,0,C(2,0,S),0) -- might be the least cardinal κ such that V$_κ$ and V$_λ$ satisfy the same Σ$_2$ statements, where λ is the least inaccessible cardinal.
C(1,0,C(2,0,S)*S,0) -- might be the least cardinal κ such that V$_κ$<$_{Σ_2}$V$_λ$, where λ is the least inaccessible cardinal.
C(C(ω,0,S),0) -- might be the proof theoretical ordinal of ZFC.
C(1,0,C(a,0,S),0) -- approximately the least cardinal κ that is not Σ$_1$(Card) definable in L$_a$ (for appropriate large *a* ) using parameters <κ. Equivalently, it is the least κ that is inaccessible in the collapse of Σ$_1$(Card) Skolem hull of L$_a$ where ordinals <κ are constants.
C(1,0,S$^+$,0) -- the least inaccessible cardinal.
C(1,0,S$^+$+S,0) -- the least κ that is a limit of κ inaccessible cardinals.
C(1,0,S$^+$*2,0) -- the least hyperinaccessible cardinal.
C(1,0,S$^+$*C(2,0,S)*S,0) -- might be the least cardinal κ such that V$_κ$<$_{Σ_2}$V$_λ$, where λ is the least Mahlo cardinal.
C(1,0,S$^+$*C(a,0,S),0) -- approximately the least cardinal κ that is Mahlo in the collapse of Σ$_1$(Card) Skolem hull of L$_a$ (for appropriate large *a*) where ordinals <κ are constants.
C(1,0,(S$^+$)$^2$,0) -- the least Mahlo cardinal.
C(1,0,(S$^+$)$^3$,0) -- the least hyper-Mahlo cardinal.
C(1,0,(S$^+$)$^{S^+}$,0) -- the least greatly Mahlo cardinal.
C(1,0,S$^{++}$,0) -- might be the least weakly compact cardinal κ (assuming representations for ordinals below κ$^{++}$ can be sufficiently robust, as otherwise the least weakly compact could be lower).
C(1,0,S$^{+(n+1)}$,0) -- might be the least Π$^n$$_1$ indescribable cardinal.
C(1,0,a,0) -- approximately the least *a*-indescribable cardinal (for appropriate large ordinal *a*).
C(1,n,0,0) -- might be the least n-subtle cardinal.

Other possibilities for n-subtle cardinals are C(n+1,0,0,0) (with C(1,a,0,0) being approximately *a*-indescribable) or C(Ω$^{n+1}$,0). While we do not see which set theoretical structures would raise the assignments this way, much uncertainty remains.

*More on cardinals after the first subtle cardinal:*
S is also the least cardinal κ such that for every S⊂κ, there is an S-indescribable cardinal <κ.
Let T = C(1,2,0,0).
C(1,1,C(1,0,a,T),0) -- (for appropriate large *a*) approximately the least *a*-indescribable limit of subtle cardinals.
C(1,1,C(1,1,0,T),0) -- the least subtle limit of subtle cardinals.

$C(1,1,C(1,1,0,T)*T^+,0)$ -- the least κ such that subtle cardinals are stationary below κ.
$C(1,1,C(1,1,0,T)*C(1,0,a,T),0)$ -- (for appropriate large *a*) approximately the least *a*-indescribable subtle cardinal.

Beyond n-subtle cardinals, we have α-subtle cardinals where instead of an n-tuple of indiscernibles, we have a well-founded tree, with the indiscernity along each branch of the tree (but not necessarily between different branches). Formally, κ is α-subtle ($α≤κ^+$) iff for every α'<2+α, club C⊂κ, and regressive $f:κ^{<ω}→κ$, there is a well-founded tree of rank α' of ordinals in C such that every branch is an increasing tuple homogeneous for f. Note that a completely ineffable cardinal is α-subtle for every α, and is below the least 1-iterable cardinal. Possible assignment (low confidence):
$C(1,α,0,0)$ -- approximately the least α-subtle cardinal (if α is not too large).
$C(1,C(2,0,α,0),0,0)$ -- approximately the least α-iterable cardinal (in a generic extension where α is countable, and assuming the cardinal is above α).
$C(1,P,0,0)$ -- (where P corresponds to a large cardinal axiom) approximately the least possible for a transitive model M of P, with every M-cardinal being a cardinal in L, with the M existing in a generic extension of L (note that M can have measurable cardinals and more).
Going further, recall that the above is in the $Σ_2$ hull of L.
$C(2,x,y,z)$ appears work similarly but with $Σ_3$ hull of L, and so on.
$C(α,0,0,0)$ -- the height of the transitive collapse of (approximately) $Σ_α$ hull of L (assuming $α<C(α,0,0,0)$).
Beyond that, using n ordinals of degree $Ω_2{}^3$, we get hulls using about n Silver indiscernibles, allowing us to approach the true $ω_1{}^L$.
$C(C(1,0,0,1,0),0,0,x)$ -- the least L-cardinal above x. At this point, we effectively have zero sharp, so successor infinite L-cardinals have effective cofinality ω (with the cofinality witnessed by the notation system).
$C(C(1,0,0,1,0)*a,0,0,0)$ -- approximately the result of iterating the sharp operation *a* times.
$Ω_2$ -- $ω_1{}^V$. We do not use $ω_1{}^{HOD}$ because we expect that (analogously to $ω_1{}^L$) in an appropriate extension of the language of set theory, there is a canonical sequence of length ω converging to $ω_1{}^{HOD}$.

For the full n=2 system, ordinals <Ω may correspond to countable ordinals, and higher ordinals correspond to countable mice with mouse order $<ε_{ω_1+1}$ (or $<ε_{Ord+1}$ among all mice). (A mouse is a generalization of $L_a$ for canonical inner models.) C(a,b) roughly corresponds to ordinal height of the least mouse of order a' (a' depends on b), where a' is obtained from *a* by representing *a* in terms of ordinals <Ω and replacing each term $a_i$ that is too large with $a_i$', where $a_i$' corresponds to a model of the description of $a_i$ with the model having definability level corresponding to the position of $a_i$ in *a*. For example, if $C(Ω^i*a,0)$ corresponds to iterating some mouse operator M *a* times, then for large *a* it may correspond to the least model of "a" that is closed under M.

Beyond n=2, as of this writing, our understanding is very limited, so we present 3 options.

The low option assumes that we are missing something, and do not reach Woodin cardinals. For example, $\mathbf{\Sigma^1_{n+3}}$ generic absoluteness is equiconsistent with n strong cardinals (with model also apparently having $\mathbf{\Sigma^1_{n+3}}$ uniformization), and it could be that some of the systems are stuck there.

For the main option, we note that in the nth system, we permit n drops, which corresponds to mice having, in some sense, complexity level n. Moreover, if ordinals correspond to real numbers, a drop corresponds to a quantifier, and mice capturing $\Sigma_n$ truth correspond to about n-2 Woodin cardinals.

An old guess (just a guess) is that:
$C(\Omega_3*a,0)$ -- enough structure for (approximately) *a*-reflective cardinals, especially for a>ω (reflective cardinals are defined in (Taranovsky 2012)).
$C(\Omega_3^{a+1},0)$ -- measurable of order *a*.

For the (aspirational) high option, we transcend Woodin cardinals, with $C(\Omega_3+\alpha,0)$ (or some other term) corresponding to a countable mouse with roughly α (or $\omega^\alpha$) Woodin cardinals. The most optimistic possibility is that ordinals of correctness n>2 correspond to (approximately) n-2-reflective cardinals (and their limits), and also that n-reflective cardinals agree with the critical sequence of an n-huge embedding. (Note that the first is about the level of expressiveness, while the second is about strength, and capturing n-reflective cardinals in V may or may not naturally correspond to reaching (approximately) n-huge cardinals analogously to how the expressivenss of second order arithmetic leads to the strength of Woodin cardinals.) The supremum of representable countable ordinals might then be $\omega_1^{HOD}$ (recall that reflective cardinals go beyond (V,∈) definability). Also, n-huge cardinals have a natural definition using a normal ultrafilter, and can also be characterized by M including a sufficient fragment of the embedding, $j''(j^{(n)}(\kappa)) \in M$.

# 7 Iteration of n-built from below

By iterating n-built from below, we get a candidate stronger notation system, but its strength and well-foundedness are unclear. The construction is similar to [8.2 Degrees of Reflection with Passthrough](#) (especially for handling the limit stages of the iteration), but with additional complications. For relevant subterms we have to set the permissible n (for n-built-from-below) so as to prevent infinite regress, which we accomplish by reading off n from the portion of term that is sufficiently stable relative to the subterm.

**Notation System:**
**Syntax:** Two constants (0, Ω) and a binary function C.
**Comparison:** Standard. (For ordinals in the standard representation written in the postfix form, the comparison is done in the lexicographical order where 'C' <

'0' < 'Ω': For example, C(C(0,0),0) < C(Ω, 0) because 000CC < 0ΩC.)
**Standard Form:**
0, Ω are standard
"C(a, b)" is standard iff
**1.** "a" and "b" are standard
**2.** b is 0 or Ω or "C(c, d)" with a≤c, and
**3.** Every a'<Ω in the representation of *a* in terms of ordinals <Ω representation of *a* is a''-n-built from below from <C(a,b). See below for definitions.
*Possible extension:* Use a-n-built from below (but n will still be based on a'').

**Definitions:**
**4.** *a* is a''-0 built from below from <b iff a<b.
   *a* is a''-n+1 built from below from <b iff *a* does not use ordinals <Ω above *a* except
      a. as a subterm of an ordinal a''-n built from below from <b or
      b. as a proper subterm of C(c,d) with c<a'' that is not in the scope of a subterm of C(c,d) that is <C(c,d).

**5.** Testing for <C(a,b) can be done in the ordinary lexicographical order.
**6.** a'' is the part of *a* that is not affected by taking the limit a'→Ω, and n is maximal such that a'' (as a string in prefix form) has a (possibly empty) prefix 'C'*i + 'C0'*n, where '+' denotes concatenation and '*' denotes repetition (essentially, n is the integer part of the least significant part of a'').
**7.** To get the limit in (6), which will call $a_{lim}$, represent *a* in prefix form, delete everything to the left of a', replace a' with Ω, add the right number of Cs on the left to make the term well-formed (not necessarily standard), and then recursively replace C(d,C(e,f)) with C(d,f) where d>e (using standard comparison; replacement order does not matter). a'' (as a string) will be the largest common suffix of a and $a_{lim}$; if a'' is empty, replace it with 0.

**Extensions:** If for representing ordinals above Ω (in terms of ordinals <Ω), we use a given C (perhaps using Ω' to go beyond CNF), the above construction applies, with the limit in (6) obtained directly using nonstandard forms above Ω (it is unclear if the computation in (7) will work for typical C).

**Variation using CNF base Ω:**
(1)-(5) are as above, except as modified below.
**8.** In constructing the subterm tree of *a* (used below), ordinals ≥Ω are represented in Cantor Normal Form base Ω, which can be accomplished using the equation C(a,b)=b+$\omega^a$ iff *a* is ≤ the least fix-point of x→$\omega^x$ above b (the inequality always holds if b≥Ω). (Also, C(a,b)<Ω iff b <Ω.)
**9.** a'' is obtained by representing *a* in CNF base Ω and deleting all terms with significance level less than or equal that of a', where $\Omega^c$*d is deleted if c or d is deleted, and '+' is deleted if d is deleted in c+d. The significance level of a term d<Ω in *a* (using CNF base Ω) is obtained by deleting everything to the right of d (but if d is inside e in $\Omega^e$*f; f is replaced by 1), changing d to Ω and subtracting *a* on the left (and converting to standard form for comparison).
**10.** n=max(m∈ℕ: ∃d∈Ord c=d+m), where c is the rightmost ordinal <Ω in CNF

base Ω representation of *a*, where c is to the left of a' and the significance level (as defined above) of c is higher than that of a' (and c is 0 if *a* has no terms of higher significance level than that of a').

**Notes:**
* In 4, the same ordinal may occur in several places in the subterm tree, and each occurrence is treated separately.
* In practice, one would use CNF base Ω for ordinals ≥Ω (also, one would likely prefer 1-variable C).
* A convenient variation is to use n+1 instead of n; it has the same strength but especially for a<Ω, more closely resembles our other systems.
* As written, if $Ω^c$*d is deleted because d is deleted, subterms of c might still be used for setting n.
* Nonstandard forms can be defined in the usual way: C(a,b) is the least ordinal >b that equals C(c,d) (standard form) for c≥a. (Or if we want to make C(a,b) total in b but partial in *a*, set C(a,C(b,c)) = C(a,c) if a>b, and after all these simplifications are complete, require the form to be standard.)
* A variation is to use a different ordinal notation system (instead of C) for ordinals ≥Ω. One can use an arbitrarily strong system as long as it has an appropriate structure for the combined system to work.
* A variation (apparently permitting fewer terms) is not to consider ordinals inside passthrough (4(b)) as candidates for a''-n built from below (4a), thus making the structure of C(c,d) above <C(c,d) irrelevant.
* *Reflection configurations:* Given our special use of the integer part, reflection configurations are only (unconditionally) defined above limit ordinals. For b+ω, configurations above <(b+ω) use a fictious limit ordinal x with b<x<b+ω. This limitation does not apply to configurations that do not use Ω (except in constants below x).
* The well-foundness of the possible extension is unclear. We still need to be strict about n so as to prevent infinite descending sequences like Ω, C(Ω+1,0), C(C(Ω+1,0)+2,0),… (each lower ordinal gets a higher permissible n).

To understand the system (using CNF base Ω), consider the segment that uses only C(a,b,c) = C(Ω*a+b,c) (with a,b,c <Ω). 4 requires that b is a-n-built from below from <C(a,b,c) where n is the integer part of a. (I use a-n in place of Ω*a-n; also in this segment *a* is automatically built from below from <C(a,b,c).) For a=n, this approximately corresponds with n built from below from c, which allows embedding of the main system (the one using $Ω_n$) into this segment with finite *a*. If the main notation system corresponds to $Z_2$+PD, for limit *a*, the level may correspond with $J_a(\mathbb{R})$, and then one repeats n-built from below construction but this time assuming that for lower levels definability has been completed, hence 4b. The sense of 4b is that C(c,d) is treated as representing C(c,d) in terms of ordinals <C(c,d), and which intermediate ordinals >C(c,d) are used is irrelevant.

The full system is a transfinite iteration of the above. Beyond $Ω^e$ for a fixed small ordinal e, CNF base Ω becomes more complex, but using significance level should be the right generalization. In requiring a' to a''-n built from below (from

<C(a,b,c)), we avoid circularity by requiring that, in a certain sense, a'' is constructed prior to a', which can be used to explain the construction of a''.

*Strength*: For all we know, if the main system corresponds to second order arithmetic or lower, the strength here might be only slightly higher than second order arithmetic (or even lower, or if above $Z_2$, just slightly below existence of $L_{\varepsilon_{\omega_1+1}}$). However, if the main notation system (that uses $\Omega_n$) reaches $Z_2$+PD, a reasonable hypothesis for the strength of this system below $\Omega^{\Omega*\omega}$ is rudimentary set theory + for every ordinal κ there are κ Woodin cardinals (with ordinals <Ω corresponding to Wadge ranks (within determinacy) that have a canonical definition in the theory). This appears to have the same strength as rudimentary set theory plus schema (n a natural number): Games on integers of countable length and projective payoff are determined, where the game consists of n rounds, with the first round having length ω and each round coding (with coding having projective complexity) the length of the next one. $\Omega^{a+1}$ may correspond with (approximately) ω*a Woodin cardinals, and $\Omega^{\Omega+a}$ with (approximately) the number of Woodin cardinals equaling the value of (ω*a)th Woodin cardinal. The full system may correspond to games (of variable countable ordinal length) with a level n round (schema given n), where level m rounds consist of level m-1 subrounds, with the number of subrounds (or number of moves for level 1 rounds) determined at the start of each round. (This strength can also be expressed using limits of Woodin cardinals, and is well below a regular limit of Woodin cardinals.)

A question regarding this assignment of strengths is that if n-built-from-below leads to around n Woodin cardinals, what happens if the system were built using (say) 5-built-from-below throughout. On one hand, it is unclear what would prevent it from accumulating Woodin cardinals (at a slower rate, getting around kα (k<5) in place of ωα Woodin cardinals at similar points), but on the other hand, the construction appears weaker than (say) 7-built-from-below (in the main system). It is possible that our conjectured strengths are wrong, with a small chance that the main notation system already reaches n-huge cardinals.

**Further extensions:**
* While the strength can be increased slightly by going beyond $\varepsilon_{\Omega+1}$, new construction principles are required to reach substantially further.
* One approach is to set n in some way that more closely approaches self-reference.
* Another approach is that in C(*a*,b,c), one would want to allow *a* to have a richer structure than b. Here is how it may look like, although the particular extension here is just an illustrative guess and might well be either ill-founded or failing to reach further. In 4 (above), for the variation using CNF base Ω, add clause
c. or the term f (a'<f<Ω) has lower significance level in e than that of a' in *a*, where e is the outermost (without intervening C) CNF base Ω superterm of f. (Thus, it is more accurate to say (significance level of a' in a)-a''-n built from below. The significance level is for CNF base Ω and is as defined above. If f is e, the level is Ω. One can also consider the restriction below (for example) $\Omega^\omega$, where the significance level of a' in $\Omega^n$a' is (essentially) n.)

**Reflection configuration version:** If iterations of n-built-from-below is weaker than expected, using 2.5 Reflection Configurations may allow the system to get the right strength. The definition combines the reflection configuration definitions for the main notation system and for Degrees of Reflection with Passthrough. Let r' denote $r_{<\Omega}$. For comparing configurations, Ω stands for λx<Ω.Ω. The relevant condition will be that for every d<Ω in the representation of *a* in terms of ordinals <Ω, r'(d) is r'(a),n-built-from-below (or r'(a''),n for the unextended version), denoted by $B_{r'(a),n}(r'(d))$.
d is a,0-built-from-below iff it is a constant (below Ω = λx<Ω.Ω).
d is a,n+1-built-from-below iff (with quantifiers ranging over $T_d$)
∀x⊑d (d < r'(x) < Ω) ∃y⊑d (r'(x)⊆r'(y)) (y<Ω ∧ $B_{a,n}$(r'(y)) ∨ parent(y)<Ω ∧ r'(y)<a).

**Future:** To get to n-huge cardinals, one would want to find how to place a corresponding embedding j in a simple framework.

# 8 Built-from-below with Passthrough

## 8.1 Definition and Properties

Here we present a strong ordinal notation system with conjectured strength beyond second order arithmetic, and below (8.2 Degrees of Reflection with Passthrough) we hope to reach beyond ZFC. The system comes in two versions:
* The basic version whose well-foundedness we prove below (using large cardinals). It behaves very similar to the extended version (which makes it exceptionally precise in imitating second order arithmetic), except that because of delayed diagonalization, we suspect that it is weak.
* The extended version that uses reflection configurations and should have the right strength.

$C_i$(a,b) will be the least ordinal of reflection level i and for that level of degree *a* above b. The sense of reflection level i is that definitions for reflection levels <i have been completed in a sense that will be formalized by the notation system. The system can be described as built-from-below with passthrough for lower reflection levels.

*Relation with the main system:* $C_i$(a,b) superficially corresponds with the main system C($\Omega_2$*i+a,b) with a,b<$\Omega_2$ and *a* and b representable without using ordinals ≥$\Omega_2^2$. However, the $C_i$ system permits additional terms, including apparently all terms with C(Ω*i+a,b) valid for the Degrees of Reflection system, such as $C_0(C_1(C_2(0,0),C_2(0,0)),0)$, and most likely, the main system misses too many terms for the above correspondence to approximately hold. In the other direction, not all valid main system terms that do not use ordinals ≥$\Omega_2^2$ become valid $C_i$ terms; an example is (courtesy of Hyp cos)
C($\Omega_2$+C(C($\Omega_2$*2+C($\Omega_2$+C($\Omega_2$,C($\Omega_2$*2,0)),0),0),C($\Omega_2$*2,0)),0).

**Syntax:** Constant 0 (the least ordinal), and function C with 3 arguments.

**Standard form:**
$C_i(a,b)$ is standard iff
  - a and b are standard
  - b is 0 or b is $C_j(c,d)$ with $(i,a) \leq (j,c)$ (lexicographically). (This is needed for uniqueness of representations, but not needed for comparison.)
  - a is built from below from $<C_i(a,b)$ with passthrough for $(C_j : j<i)$, that is *a* does not use ordinals above *a* except
    -- as a subterm of an ordinal $<C_i(a,b)$ (using standard comparison to check) or
    -- as a proper subterm of $C_j(e,f)$ where $j<i$ that is not in the scope of a subterm of $C_j(e,f)$ that is $<C_j(e,f)$. In other words, $C_j(e,f)$ is treated as representing the ordinal $C_j(e,f)$ in terms of ordinals $<C_j(e,f)$, and which intermediate ordinals $>C_j(e,f)$ are used is irrelevant.

**Comparison:** Standard where $(i,a)$ is treated as an ordinal $((i,a)<(j,b) \Leftrightarrow i<j \vee i=j \wedge a<b)$. That is (assuming standard form)
$C_i(a,b) \leq d \Rightarrow C_i(a,b) < C_j(c,d)$
$b \geq C_j(c,d) \Rightarrow C_i(a,b) > C_j(c,d)$
otherwise, $C_i(a,b) < C_j(c,d)$ iff $(i,a) < (j,c)$.

**Example:** Let $M=C_2(0,0)$, and for readability allow ordinal addition.
$C_1(M+C_1(M*2,0),0)$ is invalid because $M+C_1(M*2,0)$ uses a higher ordinal $(M*2)$.
$C_1(C_0(C_3(0,0),M),0)$ is valid because of the passthrough for $C_0$. As far as $C_1$ is concerned, we are not using $C_3(0,0)$ as an ordinal but merely as a notational shorthand to define $C_0(C_3(0,0),M)$ in terms of M.
On the other hand, in the first example, $M+C_1(M*2,0) = C_0(C_1(M*2,0),M)$ and $C_1$ "knows" what is $C_0(C_1(M*2,0),M)$ in terms of $C_1(M*2,0)$ and M but cannot handle $C_1(M*2,0)$ because $M*2$ is larger than $M+C_1(M*2,0)$.

**Proposition:** In the definition of $C_i(a,b)$, it is sufficient to allow ordinals $\leq b$ as constants. For the passthrough it does not matter whether we use $\leq f$ or $<C_j(e,f)$.
**Note:** Out of the four equivalent possibilities here, we chose the one that makes sense for reflection configurations.
**Proof:** Regarding passthrough, all terms strictly between f and $C_j(e,f)$ have degree $<(j,e)$, so passthrough applies to those terms as well. Regarding $\leq b$ vs $<C_i(a,b)$, assume contrary, and let b'>b be minimal such that *a* is built from below (with the passthrough) from $\leq b'$ (because b' is a subterm of *a*, minimality does not require proving well-foundedness). b' is $C_j(a',b'')$ with $j<i$ or $j=i \wedge a'<a$, and in both cases *a* is built from below (with the passthrough) from $<b'$, contradicting minimality.

**Building the system above an ordinal:** The system can also be built above any given ordinal α. One version is to add ordinals $<α$ as constants and to have $C_i(a,b) \geq α$. Another version with the same strength is to add ordinals $\leq α$ and (for standard forms) only use $C_i(a,b)$ for $b \geq α$; its benefit is that $C_i$ need not depend on α.

Terms that are built-from-below typically correspond with recursive terms, or to ordinals below the proof ordinal of the theory represented. Since built-from-below terms are captured by $C_0$, a natural assignment is to set $C_1(0,x)$ as the next admissible ordinal above x. Moreover such assignment is canonical in its ability to keep existing terms unchanged as one extends $C_0$ to a stronger recursive system. The relationship for higher $C_i$ and $C_{i+1}$ is analogous. We conjecture that this extensibility (properly formulated) can be used to formally define the canonical assignment of ordinals (and that the strength of the system corresponds to (or is weaker than) the assignment).

**Reflection Configuration Version**: If (as it seems likely; see proof ordinals below), the delay in the ability to do diagonalization affects the strength, our remedy is to use [2.5 Reflection Configurations](#) rather than degrees. In (for example) $C_j(c,M)$ (j≥i) as part of $C_i(a,b)$, the relevant test will be not whether c<a, but whether what c is to $C_j(c,M)$ is less than what *a* is to $C_i(a,b)$, which is expressed using reflection configurations.
**Formalization:** In the subterm tree (with terms in different positions distinguished) for C(a,b), let r(c) for c in $C_j(c,d)$ be the reflection configuration of c above $C_j(c,d)$; these are the only configurations we consider. In the definition of the system, replace "*a* does not use ordinals above *a*" with "*a* does not use reflection configurations >r(*a*)".
**Notes:**
* Out of j,e,f, only e can be bigger than $C_j(e,f)$. Thus, in the original version, it suffices to test e (for being >a), and in this form, the only difference (from using reflection configurations) is in how e is tested.
* All terms valid in the basic version are also valid here because in the relevant case, $C_i(a,b) < C_j(e,f)$.

**One variable C:** We will use [2.4 One Variable C](#), with $C'_i(\omega^{a_n}+...+\omega^{a_1}) = C_i(a_1,...,C_1(a_n,0)..))$ (standard term). Note that $C'_i(0) = 0$ and $C'_i(\omega^2) = C_i(\omega 2, 0)$.

**Conjectured Strength:** A natural conjecture is that the $C_n$ system (n>0, terms only use $C_i$ for i≤n) corresponds to $\Pi^1_n\text{-TR}_0$. (TR stands for transfinite recursion and 0 indicates limited induction.) For the conjectured natural assignment of gaps, $C_{1+i}$ returns ordinals κ such that $L_\kappa <_{\Sigma_i} \omega_1^L$ for i>0, allowing transfinite i, and $C_1$ returns admissible ordinals and their limits, and with ordinals below $C_1(0,0)$ recursive. The proof theoretical ordinal of second order arithmetic is likely $C_0(C_\omega(0,0),0) = C'_0(C'_\omega(1))$. The full system may then correspond to (canonical countable ordinals in) rudimentary set theory plus for every countable α, $L_{\omega_1+\alpha}$ exists. This has the same strength (and the same canonical countable ordinals) as second order arithmetic with comprehension for infinitary formulas, axiomatized using the truth predicate or satisfaction relation for such formulas.

**Examples (conjectured canonical assignment):**
Let $M = C_2(0,0) = C'_2(1)$, and let $x^+ = C_1(0,x)$, and $N = M^+$.

C'$_1$(a) for a<M enumerates admissible ordinals <M and their limits, C'$_1$(M+a) enumerates fix-points of the above enumeration, C'$_1$(M*a+b) corresponds to 2-variable enumeration of the fixpoints, and so on for ordinals recursive in M (that is for ordinals <M$^+$). C'$_1$(N+a) (a<N=M$^+$) enumerates recursively inaccessible ordinals <M and their limits, and so on with N being similar to Ω in [4 Degrees of Reflection](). For example,
C'$_1$(N$^{N^{N^2+N}}$*16+N$^{N^N}$*15+N$^{N}$*N*14+N$^N$*13+N$^2$*12+N*11+10) should be the 10th admissible above the 11th recursively inaccessible above the 12th recursively hyperinaccessible above the 13th recursively Mahlo above the 14th admissible limit of recursively Mahlo ordinals above the 15th Π$_3$ reflecting ordinal above the 16th Π$_3$ reflecting ordinal that is Π$_2$ reflecting onto Π$_4$-reflecting ordinals. Note that while diagonalization in "Degrees of Reflection" uses built-from-below ordinals <Ω, here we first collapse a built-from-below ordinal into an ordinal <M$^+$ (or another appropriate ordinal), and then use it (using passthrough).

Going further, for appropriate f, C'$_1$(f(M)) is roughly the least ordinal c that is f(c) stable. However, this is approximate.

* For example, C'$_1$(M$^{++}$) is the least $\Sigma^1_1$-reflecting c (and thus the least non-Gandy c) as opposed to the least c$^{++}$ stable c. Also, for non-Gandy κ, C$_1$(0,κ) appears to be the supremum of κ-recursive well-orderings of κ, with the next admissible above κ being C$_1$(0,(C$_1$(0,κ))) instead, or even higher for very strong reflection levels (or if we missed essential structure). (Also, since M is non-Gandy, N is not an admissible ordinal.)
* C'$_1$(φ$_{M^{++}}$(1)) (using Veblen φ) could be the least ordinal that is both $\Pi^1_1$ and $\Sigma^1_1$ reflecting.
* C'$_1$(φ$_{M^{+n}}$(1)) (n>2) could be the least ordinal stable up to n-1 admissibles.
* C'$_1$(M$^{+\omega}$) = C'$_1$(C$_1$(1,M)) is likely the least ordinal that for every n<ω is above an ordinal that is stable up to n admissibles.
* C'$_1$(φ$_{M^+}$(M$^{+\omega}$+1)) is likely the least ordinal that for every n is $\Pi^1_1$-reflecting onto κ$^{+n}$-stable κ . The use of φ appears analogous to the assignment of degrees for "Degrees of Reflection" above; different reflection levels can be combined and iterated similar to that assignment.
* C'$_1$(M$^{+(\omega+1)}$) is likely the least κ$^{+\omega}$-stable κ.
* For the basic version, the least ordinal stable up to a recursively inaccessible ordinal is likely C'$_1$(C'$_2$(2)$^+$). This is the because C'$_2$(2)$^+$ = C$_2$(0,M)$^+$ is required to get the first recursively inaccessible ordinal above M, I = C$_1$(C'$_2$(2)$^+$,M) (and I is the first relevant ordinal that requires C'$_2$(2)$^+$). The least ordinal $\Pi^1_1$-reflecting onto ordinals stable up to a recursively inaccessible is likely C$_1$(C'$_2$(2)$^+$+φ$_{M^+}$(I+1),0), and similarly with other ordinals. However, for the reflection configuration version, built-from-below is more permissive, and the ordinals are C'$_1$(I) and C'$_1$(φ$_{M^+}$(I+1)), respectively.
* For the basic version, the least ordinal above x<M stable up to a nonprojectible

ordinal could be $C_1(C'_2(\omega)^+,x)$, in which case the height of the least transitive model of $\Pi^1_2$-CA is $C_1(C'_2(\omega)^++1,0)$. For the reflection configuration version, use $C'_2(2)$ in place of $C'_2(\omega)$.

**Additional Intuition:** In $c = C_i(a,b)$ (for the standard choice of *a*), let $M = C_{i+1}(0,b)$. c is the least ordinal >b that satisfies (approximately) T, where T denotes the reflection properties of M up to *a* (using ordinals <c in $C_i(a,b)$ as constants), and (as a description) T works at absoluteness level i. Now, $C_i(a,b)$ is monotonic (including in *a*), and to achieve that, the standard *a* is built-from-below using functions that are sufficiently absolute between c and M. For example, $C_0(a,b)$ gives a recursive-in-b ordinal, and the ability to use passthrough for $C_0$ corresponds to $\Pi^1_1$-absoluteness for appropriate transitive models. Now, $\Sigma_i$ elementary substructures agree about $\Sigma_i$ formulas, and each increase in i corresponds to an extra quantifier, or an extra definability level for the passthrough.

**Proof ordinals in the reflection configuration version:** The use of reflection configurations can replicate the structure $\omega$ times. For example, the system with all subterms below $C'_1(2)$ corresponds to KP, but $\varepsilon_0 = C'_0(C'_1(1)) < C'_0(C_0(C'_1(2),C'_1(1))) < C'_0(C_0(C_0(C'_1(3),C'_1(2)),C'_1(1))) < ... < C'_0(C'_1(2)) =$ the proof ordinal for $\Pi^1_1$-CA$_0$, with each term in the above sequence essentially adding an extra admissible to KP. It appears that the replication happens at each admissible, leading to $C'_0(C'_1(\omega))$ being the proof ordinal for $\Pi^1_1$-CA$_0 + \omega^\omega$ admissibles.
$C'_0(N)$ is likely the proof ordinal for $\Pi^1_1$-CA$_0 + \omega$ recursively inaccessible ordinals, and similarly with $N^2$ and recursively Mahlo ordinals, and so on.
$C'_0(C'_i(2))$ is likely the proof ordinal for $\Pi^1_i$-CA$_0$.
$C'_0(C'_{i+1}(1))$ is likely the proof ordinal for $\Pi^1_i$-TR$_0$.

**Layering and delayed proof ordinals:**
   However, for the basic version, $C'_0(C'_2(n))$ for n≤$\omega$ is only the proof ordinal of $\Pi^1_1$-TR$_0$ + n recursively inaccessible ordinals. To see this, let $T_{i,c}$ (for sufficiently definable c) be the notation system where (1) no term or subterm is ≥c, (2) the *i*th recursively inaccessible ordinal is assigned to $C'_2(i)$, and (3) in place of terms <$C'_2(i)$, all ordinals <$C'_2(i)$ are used as constants. If $T_{i+1,c}$ is well-founded, then so is $T_{i,c}$ because recursiveness of $T_{i+1,c}$ combined with the recursive inaccessibility allows us to carry out an ordinal assignment. Formulated in a different and more specific way, instead of setting $C_1(0,x)$ as the next admissible ordinal, we can set it to the proof ordinal (built above x) of $\Pi^1_1$-TR$_0$ + (n-i) recursively inaccessible ordinals, where $C'_2(i)≤x<C'_2(i+1)$. Similar layering for $C'_i(n)$ also happens with higher i.
   We have not ruled out that (for appropriate a<$C'_2(\omega)$) $C_0(C'_2(\omega)+a,x)$ "unlocks"

the strength of *a* above x, but our expectation is low. $C_0(C'_2(\omega) + C'_1(n),0)$ (n>1) appears to be the proof ordinal of KP + there is a limit of recursively inaccessible ordinals and n-2 admissibles above it. A lower bound on $C_0(C_2(\omega)*\alpha,0)$ is the proof ordinal for $\Pi^1_1\text{-CA}_0$ + $\omega\alpha$ recursively inaccessible ordinals, using diagonalization to express $\alpha > C'_2(\omega)$. This lower bound (and its generalization to other limits) is insufficient to rule out correspondence with [3 Degrees of Recursive Inaccessibility](#) (with $C_0(C_{i+1}(\omega))$ being the proof ordinal of $\Pi^1_1\text{-CA}_0$ + $\omega$ recursively i-inaccessible ordinals), or between [4 Degrees of Reflection](#) and [8.2 Degrees of Reflection with Passthrough](#).

Ordinal notation systems can be characterized by (1) the kind of diagonalization they permit, and (2) by how high up you have to go to carry out diagonalization. On the first point, we have enough structure for second order arithmetic. However, on the second point, the notation system (without reflection configurations) resembles [3 Degrees of Recursive Inaccessibility](#). To give an example, while getting a diagonalization of admissibles (in typical systems) can be done above the next recursively inaccessible ordinal, here we have to wait until $C_2(0,0)$. Note that even if $C_0$ is too weak in the full system, the canonicity of the above (partially given) assignment remains (besides correspondence with the reflection configuration version) in the ability to keep the assignment unchanged as one extends $C_0$, including in [8.2 Degrees of Reflection with Passthrough](#).

## 8.2 Degrees of Reflection with Passthrough

The built-from-below with passthrough for lower levels system can be extended by adding more variables. One extension (not the one below) which might be too restrictive is as follows. With $C_{a,b}(c,d)$, c may use passthrough for $C_{a',b'}$ where (a',b')<(a,b) (lexicographically) and b may use passthrough for $C_{a',b'}$ where a'<a. To go further, combine a, b, and c into an ordinal using a large ordinal $\Omega$ and naturally (using CNF base $\Omega$) extend the combination up to $\varepsilon_{\Omega+1}$. The result will be the "Iteration of n-built from below" notation system with "n-built from below" replaced with "built from below". However, the above is probably too restrictive, and in $C_{a,b}(c,d)$ we should allow passthrough for b where (a',b')<(a,b). A natural extension up to $\varepsilon_{\Omega+1}$ (or another ordinal) is the following.

**Degrees of Reflection with Passthrough:**
*Syntax:* 0 (the least ordinal), C (partial binary function), $\Omega$ (a large ordinal), **O** (ordinal notation system for ordinals above $\Omega$ in terms of ordinals <$\Omega$).
*Choice of O:* CNF base $\Omega$; another choice is C; any system that is well-founded above a generic ordinal should work.
*Description:* C is as in General Notation. C(a,b)<$\Omega$ when b<$\Omega$. In C(a,b), each d in the representation of *a* in terms of ordinals <$\Omega$ is built from below from <C(a,b) allowing passthrough for C(a',b') where a'<a. This means that d does not use ordinals x with d < x < $\Omega$ except in the scope of an ordinal <C(a,b) or as permitted by passthrough. Passthrough means that C(a',b') is treated as representing C(a',b') in terms of ordinals ≤b', and the intermediate ordinals used are not checked for being <d.

**Notes:**
* Like in the $C_i(a,b)$ system, in built-from-below from <C(a,b), we can replace "<C(a,b)" with "≤b". Similarly, for the passthrough, it does not matter whether we use ≤b' or <C(a',b').
* The passthrough is the only difference from [4 Degrees of Reflection](). The only difference in Comparison Algorithm is that the condition using $T_a$ is now:
∀x∈$T_a$ ∀y⊑x (x<y<Ω) ∃z⊒y (z<Ω) (z⊐x ∨ z < C(a, b) ∨ z⊐y ∧ d(z)<a ∧ ∀t (y⊑t⊑z) t≥z)
where d(C(e,f))=e (for standard C(e,f)) and d(z)=∞ otherwise. "z⊐y ∧ d(z)<a ∧ ∀t (y⊑t⊑z) t≥z" is the passthrough condition. I did not verify conversion of nonstandard forms.
* Omitting z⊐y (thus permitting z=y in passthrough) would make the system ill-founded. An example infinite descending sequence is (courtesy of Hyp cos) $a_1$ = M = C(Ω*2,0) and $a_{n+1}$ = C(Ω+M*$a_n$,0).
* While it permits many additional terms, we have not ruled out that (for typical **O**), it has the same strength as Degrees of Reflection. If that is the case, see the reflection configuration version below for a much stronger system (whose well-foundedness is unclear).

If the system is well-founded (which is not yet proved), a natural conjecture is that it is similar to and has the same strength as the main system for n=2, or a hypothetical passthrough variation of the main system for n=2. See [6 Beyond Second Arithmetic]() (above) for analysis.

An example is C(a,b,c,d) with (a,b,c) compared lexicographically. Note that allowing passthrough for b using (a,b') (b'<b) allows b to use large values of c' even when using *a*. While seemingly unrestricted, this only appears to be useful in the presence of large enough gaps that there is no sign of infinite regress. The passthrough for c using (a,b,c') (c'<c) is only relevant when b>c (note that a < C(a,b,c,d) here) given the built-from-below condition.

**Proposition:** The $C_i(a,b)$ system above permits the same terms that are <Ω (and with the same comparison) as degrees of reflection with passthrough using pairing for O and identifying $C_i(a,b)$ with C(Ω*i+a,b).
**Proof:** The only difference not addressed above is passthrough for (i',a') with i'=i and a'<a, but this is vacuous. Assume contrary, and consider a shortest length counterexample C(Ω*i+a,b) with a prohibited subterm c' = $C_i$(a',b') < a. Because c' is valid in the $C_i(a,b)$ system and c'<a, the inability to use passthrough for c' does not invalidate the term, contradicting the assumption.
*Note:* The proposition relies on i < max(a,C(Ω*i+a,b)) and does not fully apply to extensions. A counterexample is C(Ω*i+C(Ω*i,j),0) with i = C($Ω^2$*2,0) and j = C($Ω^2$,0).

*Restricted passthrough:* To formalize a framework for restrictions, in addition to **O**, one specifies a restriction p on passthrough: Given *a* and a subterm (including position) d<Ω in the O representation of *a*, return p(a,d). The condition using ∀x∈$T_a$ would use d(z)<p(a,x) in place of d(z)<a. For example, the multivariable built-from-below with passthrough for lower levels would use p($Ω^i*a_i$+...+$Ω*a_1$+$a_0$, (j,$a_j$)) = $Ω^i*a_i$+...+$Ω^{j+1}*a_{j+1}$ (this is 0 if i=j; in p(..), j is used to indicate the position of $a_j$). If the unrestricted system (i.e. p(a,..)=a) were ill-founded (or just inconveniently strong), then between a safe choice (used above and in Iteration of n-built-from-below) and the ill-foundedness, there would be a spectrum of possibilities, with degrees of reflection with passthrough acting as a framework rather than a single system.

*Equivalent definition of Degrees of Reflection with Passthrough:*
\* The subterm tree for a term d < Ω in terms of ordinals <c with passthrough for <a (and treating ordinals ≥Ω as syntactic constructs) is as follows.
d is the root node.
Given d' in the tree:
d'<c (or d'=0) -- d' is the leaf node.
Otherwise, d' = C(a',b') (standard form), and
- if a'<a, the immediate subterms are the terms <d' (equivalently, leaf terms <Ω) in the standard representation of d' in terms of ordinals <d'. This is the passthrough condition.
- else the immediate subterms are b' and the terms <Ω in the standard representation of *a*' in terms of ordinals <Ω (note that the representation is trivial if a'<Ω).
\* C(a,b) is standard if the standard conditions are met (a and b are standard, and b is 0 or Ω (if using C above Ω) or C(c,d) with a≤c) and for every d<Ω in the representation of *a* in terms of ordinals <Ω, for every term f in the subterm tree for d in terms of ordinals <C(a,b) with passthrough for <a, f≤d.
*Note:* Restricted passthrough is defined analogously, using "... with passthrough for <p(a,d) ...". (Also, the use of <C(a,b) allows even a very restricted passthrough such as p(a,..)=0.)

**Reflection Configuration Version:** Using [2.5 Reflection Configurations](), the reflection configuration version of the system is obtained by replacing the condition using $T_a$ with:
∀x∈$T_a$ ∀y⊑x ($r_{<Ω}$(x) < $r_{<Ω}$(y) ∧ y<Ω) ∃z⊒y (z<Ω) (z⊐x ∨ z < C(a, b) ∨ z⊐y ∧ $r_{<Ω}$(d(z)) < $r_{<Ω}$(a) ∧ ∀t (y⊑t⊑z) t≥z)
(the only change is the use of $r_{<Ω}$).
\* Equivalently, the above definition using the passthrough subterm tree is modified as follows.
- Replace a'<a with r(a',d') < r(a,C(a,b)); the subterm tree is otherwise unchanged.
- Replace f≤d with r(f,par(f)) ≤ r(d,C(a,b)), where par(f) is the parent of f in the above passthrough subterm tree; the parent of d is C(a,b).
*Note:* Subterms in different positions are distinguished. Use standard comparison for C(a,b) (as if nonstandard C does not exist).

## 8.3  Proof of Well-foundedness

We first present a proof in ZFC of well-foundedness of $C_0,C_1,C_2$ system, and then (using a different approach) a proof of well-foundedness of the full $C_i$ system using a measurable cardinal. Then, we will build on that proof to prove well-foundedness of an extension. The proofs are for the version that does not use reflection configurations.

## $C_0$, $C_1$, $C_2$ system

**Theorem:** $C_0,C_1,C_2$ system is well-founded, even if built above an arbitrary ordinal.
**Proof:**
Work in the constructible universe L.
Set $C_2(a,b)=\omega_{b'+\omega^a}$ where b' is the maximal ordinal such that $\omega_{b'} \leq b$ (and b'=0 if $b<\omega_1$).
(Note: Instead of $C_2$, we could have used a function f to enumerate uncountable cardinals (in L). For the notation system, if 'a' is standard, then so is 'f(a)'; $f(a)<f(b) \Leftrightarrow a<b$, $f(a)<C_0(c,d) \Leftrightarrow f(a) \leq d$, $f(a)<C_1(c,d) \Leftrightarrow f(a) \leq d$.)
Let $C_1(a,b)$ be the least c>b that is not $\Sigma_1$(Card) definable in $L_{\kappa+a}$, allowing ordinals <c as constants, where κ is the least uncountable cardinal >b and Card is the predicate for cardinals (in L).
$C_0(a,b)$, when defined, is the least c>b such that a'<a ∧ b'<c ∧ ($C_0(a',b')$ is defined) $\Rightarrow C_0(a',b')<c$.
$C_0(a,b)$ is defined iff *a* can be represented using $C_0$, $C_1$, $C_2$ and ordinals $<C_0(a,b)$ as constants such that
  - all intermediate ordinals are below a.
  - each use of $C_0$ and $C_1$ is standard in the following sense:
    -- if $C_0(c,d)$ is used as a subterm, then c does not use ordinals ≥c except in the scope of an ordinal $<C_0(c,d)$
    -- if $C_1(c,d)$ is used as a subterm, then c does not use ordinals ≥c except in the scope of an ordinal $<C_1(c,d)$, or as a subterm of $C_0(e,f)$ that is not in the scope of a subterm g≤f of $C_0(e,f)$.

*Proposition:* If *a* is below the least fix-point of $x \to \omega^x$, then $C_0(a,b)=b+\omega^a$.
*Proof:* Simple, by recursion on a.
*Proposition:* For every a, there is b' such that $C_0(a,b)$ is defined iff b≥b'.
*Proof:* Choose the least b' such that *a* is representable (in the required sense) using ordinals $<C_0(a,b')$.
*Proposition:* $C_0(a,b)$ is monotonic in b, and for the class of *a* such that $C_0(a,b)$ is defined, strictly monotonic and continuous in a. When defined, $C_0(a,b)$ is the least point in the range of $x \to C_0(a,x)$ above b.
*Proof:* By construction of $C_0$.
*Proposition:* Standard comparison for $C_0$: $C_0(a,b) < C_0(c,d)$ iff $C_0(a,b) \leq d$ or (b < $C_0(c,d)$ and a<c) (assuming both are defined).
*Proof:* Follows from the above properties.

*Proposition:* In the definition of $C_0$, one could replace $<C_0(a,b)$ with $\leq b$.
*Proof:* By recursion on a. An ordinal c with $b<c<C_0(a,b)$ equals $C_0(a_1,C_0(a_2,...C_0(a_n,b)...))$ with each $a_i<a$.

*Proposition:* $C_1(a,b)$ is monotonic in *a* and b, and continuous in *a*.
*Proof:* By construction of $C_1$; monotonicity is the reason we used $\Sigma_1$ definability. Continuity holds because for each *a* the class of ordinals that are not sufficiently definable is closed except at cardinals.
*Definition:* a is maximal in $C_1(a,b)$ iff $C_1(a+1,b)>C_1(a,b)$.
*Proposition:* *a* is maximal in $C_1(a,b)$ iff *a* is $\Sigma_1$ definable in $L_{K+a+1}$, allowing ordinals $<C(a,b)$ and cardinals as constants.
*Proof:* Given a $\Sigma_1$-definition of *a*, the definition of $C_1(a,b)$ becomes $\Sigma_1$ in $L_{K+a+1}$, hence *a* is maximal. Conversely, if *a* is maximal, then $C_1(a,b)$ is by definition $\Sigma_1$ definable in $L_{K+a+1}$, and this definition can be converted into a $\Sigma_1$-definition of a.

These properties lead to the standard comparison for the notation system provided that we can show two things:
* $C_1(0,b)>C_0(a,b)$
* *a* is maximal in $C_1(a,b)$ when $C_1(a,b)$ is standard.

Suppose that we can show that $C_0(c,d)$ is recursive (or just $\Sigma_1$) in field(c:d) where field(c:d) consists of d and ordinals $\leq d$ obtained by parsing term c into ordinals $\leq d$. Then the above conditions are met since $C_1(a,b)$ is above all recursive-in-b ordinals, and *a* has a $\Sigma_1$ definition obtained by parsing its notation and using the values of $C_0(c,d)$.

*Proposition:* Standard comparison works in the system.
*Proof:* Suppose that standard comparison works below a, that is for every base ordinal b<a, when all subterms are <a, the comparison works. Allow terms using a. $C_0(c,d)$ equals d plus the order type of the system between d and $C_0(c,d)$. The built-from-below condition ensures that all ordinals $<C_0(c,d)$ are representable in terms of ordinals $\leq d$ while using only ordinals <a, thus the comparison is standard and hence the system is recursive in field(c:d) as required.
**End of Proof**

**$C_i$ system**

**Theorem:** Assuming an inner model with a measurable cardinal, the $C_i$ ordinal notation system is well-founded, even if built above an arbitrary ordinal.
**Proof:**
Let K be the core model below a measurable cardinal, and let $I_\alpha$ be the predicate for $\omega(\alpha+1)$-reflective cardinals for K. See (Taranovsky 2012) (including section Iterations of Reflectiveness) for background, but the properties we will use are the following:
* Elements (and finite increasing n-tuples of elements) of $I_\alpha$ are good

indiscernibles in $(K, \in, \alpha, (I_\beta)_{\beta<\alpha})$: $\varphi(x) \Leftrightarrow \varphi(y)$ provided that x and y are finite subsets of $I_\alpha$ of the same cardinality, where $\varphi$ is expressible in the model, and allowing sets in $V_{\min(x \cup y)}{}^K$ as parameters. Access to the restricted I is allowed as a two variable function (so one can quantify over $\beta$ in $I_\beta$).
* $\alpha<\beta \Rightarrow I_\beta \subset I_\alpha$
* K satisfies separation and replacement for predicates that use I.

Let $F_\alpha$ consist of functions that are definable in $(K, \in, (I_\beta)_{\beta<\alpha})$ (allowing ordinals $\leq \alpha$ as parameters, I is treated as a two variable function), with each element of $F_\alpha$ taking a tuple of ordinals as input and returning an ordinal at least as large as the supremum of tuple.

When defined, $C_i(a,b)$ is the least ordinal c in $I_i$ such that $C_i(a',b')<c$ if $a'<a$ and $b'<c$ and $C_i(a',b')$ is defined, and
$C_i(a,b)$ is defined iff *a* is built from below from $<C_i(a,b)$ in the ith extension of the notation system. That is *a* must be representable as a term in that system such that all subterms are $\leq a$ but allowing ordinals $<C_i(a,b)$ as constants. In representing *a*, the ith extension uses functions in $F_i$ and C with $j\geq i$ in $C_j(b,c)$ where b does not use ordinals $>b$ except in a term inside $C_l(g,h)$ ($l<j$) that is not in a subterm of $C_l(g,h)$ that is $\leq h$, where a subterm includes its position in the term. (By recursion on a, this definition is noncontradictory and complete, except that we have not proved that ith extension is a really an extension.)

*Proposition:*
1. $C_i(a,b)$ is the least ordinal in the range of $x \to C_i(a,x)$ above b.
2. $C_i(a,b)$ is strictly monotonic in *a* (when it is defined).
3. For every i and a, there is minimum b' such that $C_i(a,b)$ is defined iff $b \geq b'$
4. If $i<j$ and $b < C_j(c,d)$, then $C_i(a,b) < C_j(c,d)$
5. $C_i(a,b) < C_j(c,d)$ iff $C_i(a,b) \leq d \lor b \leq C_j(c,d) \land (i,a)<(j,c)$
*Proof:* 1-3 standard. For 4, the condition $C_i(a',b')<c$ corresponds to a closed unbounded set of possible c, and elements of $I_i$ are stationary below every element of $I_j$. Statement 5 follows from the above properties and the definition.

Thus, the $C_i$ ordinal notation system is well-founded provided that we can show that standard terms are valid in the constructed system. The only problem is that $C_i(a,b)$ may use $C_j(c,d)<a$ where $j<i$ and c uses ordinals $>a$ in a restricted way. The solution is to compute $C_j(c,d)$ without computing term c as an ordinal.

Let d'' be the maximal (as an ordinal or under standard comparison) subterm of $C_j(c,d)$ that is $<C_j(c,d)$; d'' is expressible using ordinals $\leq d$ and $C_{j'}$ for $j' \leq j$. Replace each $C_k(e,f)>d''$ with $C_k(e,d'')$ where $f<d''$, and represent $C(c,d)$ in terms of subterms $\leq d''$, so the representation will be a code in the standard form for the notation above d'', plus a list of ordinals $\leq d''$. We want to show that for some function in $F_i$, its application to the code and the list ordinals equals $C_j(c,d)$ (which will complete the proof).

Recursive test for whether $e = C_j(c,d)$:
* $\forall e'<e\ C_j(c,d) \neq e'$, and
* For every $c'<c$ and $d \leq d'<e$, $C_j(c',d')<e$ where c' uses jth extension of the system.

By recursion, this solves the problem provided that we can test whether c'<c and $C_j(c',d')$ is valid, and we can be overstrict with validity provided that at least one representation passes the test.

Overstrict validity test of $C_j(c',d')$ (ordinals $\leq d'$ are constants):
  For every subterm $C_k(e,f)$:
    - $f \geq d'$
    - if $f<C_j(c',d')$, then $f=d'$
    - $e \leq c'$
    - e does not use ordinals above e except in a term inside $C_l(g,h)$ (l<k) that is not in a subterm of $C_l(g,h)$ that is $\leq h$, where a subterm includes its position in the term.

Comparison when doing the validity test:
  - two C terms - standard
  - C term and an ordinal (ordinals $\leq d'$ and their values after functions in $F_j$) -- C term is greater
  - otherwise, represent both terms in terms of C terms (do not parse the C terms yet) and ordinals. Pairwise compare all the C-terms. Assign the C-terms to arbitrary ordinals in $I_j$ >d' consistent with the comparison. Compute both ordinals and compare. This is the point where we use the indiscernibility of $I_j$.

To compare c' with c (assuming $C_j(c',d')$ and $C_j(c,d)$ were tested valid and $d \leq d'<C(c,d)$):
- Use the comparison above. Even though d' need not equal d, the comparison of C terms with ordinals works since c does not have C-terms >d that are $\leq d'$ (if it did, we would have replaced d with d'' above).

The procedure could fail in three ways:
  - We got comparison as ordinals of terms valid for $C_j$ system wrong -- not possible given the analysis above.
  - We failed to include a valid ordinal -- not possible since (apart from a strictness in representations), we fully covered the definition of $C_j$.
  - We included an invalid ordinal.

Suppose that the computation fails (is wrong) for $C_i$. Then there is invalid $C_i(a,b)$ in the ith extension which we treated as valid. Let $C_j(c,d)$ be an invalid subterm of $C_i(a,b)$ with j,c,d valid. Then j>i, and c has a subterm $C_k(e,f)$ ($i \leq k<j$) for which our computation fails. This means that there is invalid $C_k(e',f')$ that we treated as valid with $e' \leq e<a$ (and $f' \geq f$). Iterating, we get an infinite descending sequence of ordinals (since e and, for each iteration, its analogue are valid), which completes the proof.

*Note:* If we can rule out overflows (where $C_i(a,b) \geq C_{i+1}(0,b)$), the measurability assumption can be weakened to "starting at zero, the sharp operation can be iterated any ordinal number of times".

**Reducing the assumption:** $\Sigma_i$ elementary substructures of $L_{\omega_1^L}$ lack indiscernibility but it may be possible to restrict $F_i$ to a class of functions that guarantee the invariance of comparison. To check validity and compare $C_i(a,b)$ with $C_i(a',b)$ (a and a' are given as codes, with each $C_{j(c,d)}$ having $d \geq b$ and with ordinals $\leq b$ used as constants), it is sufficient to know the codes and the relative ordering of the ordinals $\leq b$ that we used, and it may be possible to restrict $F_i$ to an analogous class of functions (though perhaps with special conditions on b or other ordinals). One would also need sufficient definability of $C_i$ to prove $C_i(a,b) < C_{i+1}(0,b)$. Intuitively, validity of $C_i(a,b)$ implies that *a* has a canonical definition of a certain type, which we can then 'drop' in an order preserving way to get a similar structure at a lower level.

**A fragment of Degrees of Reflection with Passthrough**

We only have a proof of well-foundedness of (a variation on) a small fragment of the system.

**Theorem:** Assume that there is an inner model $1+\alpha$ measurable cardinals, where $\alpha$ is an ordinal. The $C_i(a,b)$ system is well-founded even if augmented with $\alpha$ constants with (1) there is no constant c (c>0) with $b<c \leq C_i(a,b)$ and (2) in $C_i(a,b)$, i is built from below from $<C_i(a,b)$ with passthrough for lower values of i.

**Corollary:** Suppose that for every $\kappa < \min(\delta: L_\delta \prec_{\Sigma_1} L)$ there is a transitive model of "ZFC + $\kappa$ measurable cardinals". Then for every $\alpha$, the $C_i(a,b)$ system with $\alpha$ constants is well-founded.

**Proof of the Corollary:** This follows from $\Sigma^1_2$ correctness of $L_\delta$.

**Notes:**
* This theorem/proof is only for the basic (rather than reflection configuration) version, so (despite the impression below) the system could be weak (possibly at the level of many recursively Mahlo ordinals).
* One possibility is that under conjectured canonical assignment, each constant corresponds to an uncountable cardinal in L. The notes below assume that to be the case. However, as of this writing, we cannot rule out a lower assignment, including one where each constant simply corresponds to an ordinal not $\Sigma_1(\omega_1^L)$ definable in $L_{(\omega_1^L * 2)}$ from lower ordinals. However, the latter seems too low (unless we fail to reach $Z_2$) as (using $\omega_1^L$ for the first constant >0), $C(\omega_1^L,b,c)$ (with $c<\omega_1^L$) appears to correspond to $\kappa$ with $L_\kappa \prec_{\Sigma_\kappa} \omega_1^L$ (and thus not $\Sigma_1(\omega_1^L)$ definable in $L_{(\omega_1^L + \kappa)}$), and thus appear sufficient to capture ordinals $\Sigma_1(\omega_1^L)$

definable in $L_{(\omega_1^L*2)}$.

* Because of layering between different constants, for finite α, the strength is only expected to be around existence of α ordinals that are uncountable up to an admissible ordinal. I do not know if a similar collapse of strength also happens for infinite α.

* Unlike in the countable case, the least ordinal c such that $L_c$ and $L_{\omega_2^L}$ satisfy the same $\Sigma_1$(Card) statements is below the least $\Sigma_1$(Card) elementary substructure of $L_{\omega_2^L}$, and we have not determined the exact canonical assignment.

* The theorem and the proof suggest that in the table of cardinals and their recursive analogs, regular uncountable cardinals in L are recursive analogs of measurable cardinals.

* By iterating the construction (α → order type of the system) ω times (starting at 0), the ordinal we get $(C(C(\Omega^2+C(\Omega,0),0),0))$ is probably the proof theoretical ordinal for rudimentary set theory plus $\forall \alpha < \omega_1^{CK} \ \exists s = \omega_\alpha^L$. This is stronger than Zermelo set theory but weaker than Borel determinacy. Iterating it another time (for a total of ω+1 times) is equivalent to using a function to enumerate all $C_1(0,0)$ constants and should give us all canonical ordinals in the theory.

**Proof of the Theorem:**
The proof is analogous to the proof of well-foundedness of ordinary $C_i(a,b)$, but there are key differences which we will now cover in detail. Work in M, where M is the minimal inner model with 1+α measurable cardinals. The constants will be measurable cardinals and their limits, excluding the supremum of the measurable cardinals. Let $\kappa_m$ be the second cardinal above the supremum of measurable cardinals. The exact fine structure is unimportant, but we can for example set $M_a$ to $L_a[U]$ where U codes the sequence of ultrafilters.

*Definition:* $I_i$ (below the supremum of measurable cardinals) is defined as follows: $\lambda \in I_i$ iff λ and κ satisfy the same $\Sigma_1(\{j(M_{\kappa_m})\})$ predicates in $j(M_{\kappa_m+i+1})$ with parameters in $M_\lambda$, where κ is the least measurable cardinal (in M) above λ and j is the corresponding embedding.

We will use only the following properties of $I_i$.
*Definition:* getmin(i) is the least ordinal such that i is $\Sigma_1(\{M_{\kappa_m}\}, \text{getmin}(i))$ definable in $M_{\kappa_m+i+1}$.
*Proposition:*
* For all measurable κ, $I_i$ has measure 1 below κ.
* $I_i$ is uniformly-in-i $\Sigma_1(\{M_{\kappa_m}\}, i)$ definable in $M_{\kappa_m+i+2}$.
* $i<j \Rightarrow I_j\setminus(\text{getmin}(i)+1) \subset I_i$
* For all finite subsets s and t of $I_i$ with ∀(measurable κ) |s∩κ|=|t∩κ| and $\Sigma_1(\{M_{\kappa_m}\})$ φ with two free variables, $\forall x \in M_{\min(s\cup t)} \ M_{\kappa_m+i+1} \vDash (\varphi(x,s) \Leftrightarrow \varphi(x,t))$.
*Proof:*
The indiscernity holds because for every $\lambda \in I_i$ below measurable κ (κ is the least measurable above λ) for every finite increasing iterate M' of M with the first embedding using κ, (κ,s) agrees with (λ,s) for the relevant class of formulas where s is the set of critical points >κ used in the embeddings. Here, increasing means that the critical points are increasing. The agreement holds because the embeddings are sufficiently definable in $M_{\kappa_m+i+1}$.

*Notes:*
* Elements of $I_i$ are effectively indiscernibles except as to the number of measurable cardinals below the element.
* Our setup ensures that above getmin(i) elements of $I_j$ are indiscernibles even for predicates that use $I_i$ for i<j.

The ith extension is defined analogously to the proof of the previous theorem. Both i and *a* need to be built from below from <$C_i(a,b)$. Instead of j,a,b→$C_j(a,b)$ for j<i, we will allow x→y=F(n,x) when:
x is a tuple of ordinals.
y is $\Sigma_1(\{M_{\kappa_m}\},x)$ definable in $M_{\kappa_m+i+1}$; n is a code for the definition.
y is an ordinal >max(x).
there is no measurable κ with max(x)<κ≤y.

An equivalent condition on i is that getmin(i)<$C_i(a,b)$. This is easily seen equivalent to i being built from from below from <$C_i(a,b)$, allowing F as above.

As before, to show that standard terms are valid in the system we constructed, we propose computation of $C_i(a,b)$ where *a* is given as a term (in the ith extension) rather than ordinal. The comparison of terms is essentially as before. Since getmin(i)<$C_i(a,b)$, comparison with i works. The comparison between a C or F term and one of the measurable constants is trivial. When assigning C-terms to indiscernibles, each indiscernible must have the right number of measurable cardinals below it.

Suppose that the computation fails (is wrong) for $C_i$. Then there is invalid $C_i(a,b)$ in the ith extension which (in our computation) we treated as valid. Let $C_j(c,d)$ be an invalid subterm of $C_i(a,b)$ with j,c,d valid. Then, j or c has a subterm $C_k(e,f)$ (i≤k<j) for which our computation fails. (Otherwise, $C_k(e,f)$ would be sufficiently definable. Also, if k≥j, the use of $C_k(e,f)$ would not prevent $C_j(c,d)$ from being valid even if our computation of $C_k(e,f)$ fails.) This means that there is invalid $C_k(e',f')$ that we treated as valid with e'≤e<a (and f'≥f). Iterating, we get an infinite descending sequence of ordinals (since e and, for each iteration, its analogue are valid), which completes the proof.

## 8.4  Using Canonical Definitions

Here is a technique that may permit getting well-foundedness of C in the weakest system possible, though as of this writing, we do not have a proof of correctness, and it is unclear whether the definitions give us the desired C. The key idea is that instead of defining C using complex ordinal recursion, we define C directly using canonical definitions of its arguments, with the twist that the arguments are interpreted not just in V but in other appropriate models. The systems can be built above an arbitrary ordinal, but otherwise, only standard terms are valid. The assignments below leave more gaps than is canonical; the canonical assignment remains unclear.

*Typical General Construction:*
Ordinals $\geq \Omega$ are treated as syntactic constructs.
A valid term $C(a,b)<\Omega$ is evaluated as the least ordinal $\kappa>b$ such that there is a structure M, and an ordinal (or well-ordering) a':
* M is based on $\kappa$ (each notion makes this specific) and satisfies a correctness/closure condition based on a'.
* a' is the result of
  - taking the definition of *a* and (outside of M) replacing terms $\leq b$ with ordinals. $C(c,d)$ for $d<b$ is replaced with its value, or if $C(c,d)>b$, with $C(c,b)$ (below, the later is done implicitly).
  - applying the resulting definition inside M. (If the application fails, this M is ruled out.)

By induction on the complexity of terms, the definition (once specified) is complete. Note that applying the definition inside M may lead to a computation in M'$\in$M (with appropriate terms precomputed in M), and so on.

We start with the three-variable C, using $C(i,a,b)$ in place of $C_i(a,b)$. We can optionally require that $C(0,a,b) = b+\omega^a$ when 'a' is below the least fix-point $x \to \omega^x$ above b, but this appears unnecessary to get the right strength. We allow the system to be built above an arbitrary ordinal (within bounds), and similarly for the four variable and other extensions.

Here is one choice of assignment. For a valid term, $C(a,b,c)$ will be the lowest ordinal $\kappa>\max(a,c)$ such that $L_\kappa$ is a $\Sigma_{1+a}$ elementary substructure of $L_{\omega_1^L}$ and some level of L satisfies "$\kappa$ is $\omega_1$ and $\kappa+a+b$ exists" where "b" is the canonical definition of b above c. The canonical definition is the one using C and ordinals $\leq c$ as constants (see Typical General Construction above).

A more optimal choice that uses minimal assumptions for the fragment corresponding to $Z_2$ is the following: $C(a,b,c)$ is the least $a*\kappa+b'$-stable $\kappa>\max(a,c)$ with $L_\kappa$ a $\Sigma_a$ elementary substructure of $L_{\omega_1^L}$ where b' is b computed inside $L_\kappa$ but with subterms $\leq c$ replaced by constants (outside of $L_\kappa$). A further step to optimality is to use $\Sigma_{a-1}$ elementarity for finite *a*, using admissible ordinals for a=1, and using no gaps (if all ordinals $\leq c$ are permitted) for $C(0,b,c)$. To complete the optimality, one would apply a construction similar to [4.3 Assignment of Degrees](#), but it is unclear how to extend it to the high strength level here.

*Degrees of reflection with passthrough:* While the system is likely stronger than this (unless it is below $Z_2$), if it were defined more narrowly so as to correspond to $L_{\varepsilon_{\omega_1^L+1}}$, the following construction might work: $C(a,b)$ is the least $\kappa$ with $L_{\kappa+a'}$ elementarily embeddable into some $L_{\omega_1^L+a''}$ with critical point $\kappa$, where a' is obtained from *a* by representing *a* in terms of ordinals $<\Omega$, replacing terms $\leq b$ with constants, and then interpreting the resulting expression in $L_K$ (instead of $L_{\omega_1}{}^L$). Also, the assignment of individual ordinals can be optimized by separating $a$ as $a_1*\Omega+b_1$, using $\Sigma_{a_1}$ elementary substructures (using diagonalization to get $a_1 \geq \omega_1$), and using appropriate levels of stability for $b_1$. It is unclear whether this interpretation satisfies $C(\Omega*C(\Omega^2,0)+C(\Omega^3,0),0) < C(\Omega*C(\Omega^2+1,0),0)$ (a typical example in Degrees of Reflection with Passthrough). If it does and there are no other issues, that would suggest that the strength of Degrees of Reflection with Passthrough is only rudimentary set theory + $\forall \alpha$ "$L_{\omega_1}\alpha$ exists".

**Four Variable C**

For the four variable extension $C(a,b,c,d)$ (based on degrees of reflection with passthrough), the following might work. Even if the above construction works, the ideas here might be applicable to other systems. Roughly, $a$ corresponds to the level of indiscernibles preserved, b to the correctness for that level, and c to the object being modeled. Let $\Omega$ be a sufficiently-closed ordinal, and let $I_a$ be indiscernibles of level $a<\Omega$ for $V_\Omega$; for each $a$, $\Omega$ is the maximum element of $I_a$. $\forall b>a\ I_b \subset I_a$. Furthermore, $I_a$ are good indiscernibles for $V_\Omega$ augmented with (as a single symbol) $(I_{a'}:a'<a)$. Here, 'good' means that the indiscernibility holds even if allow parameters whose rank (as a set) is below the least indiscernible in a tuple. Regarding $\Omega$, require that for every $a<\Omega$ such that $V_{\Omega+a}$ exists (this qualification is relevant for weak theories), $\Omega$ is a-I-indescribable: For every $U \subset V_\Omega$, there is $\Omega'<\Omega$ such that $(V_{\Omega+a},\in,I,U)$ and $(V_{\Omega'+a},\in,I,U)$ satisfy the same statements with parameters in $V_{\Omega'}$. Since $\forall b<\Omega\ \Omega \in I_b$, I-indescribability (for $a>0$) ensures that $\forall b<\Omega'\ \Omega' \in I_b$.

For a valid term, $C(a,b,c,d)$ will be the least $\Omega'>d$ such that
- $V_{\Omega'+b'}$ and $V_{\Omega+b'}$ satisfy the same statements with parameters in $V_{\Omega'}$, allowing $(I_{a'}:a'\leq a)$ as a (single) predicate symbol. Here b' is the least ordinal such that for some I' agreeing with $I_{a'}$ for $a'\leq a$, $(\Omega,I')$ satisfies the definition of $(\Omega,I)$ in $V_{\Omega+b'+1}$, with b' obtained by taking the C-definition of b above d (that is using ordinals $\leq d$ as constants), and replacing I with I'.
- For some I' agreeing with $I_{a'}$ for $a'\leq a$, $(\Omega',I')$ satisfies the definition of $(\Omega,I)$ in $V_{\Omega'+b'+c'+1}$. Here, c' is given by the C-definition of c above d, using $(\Omega',I')$ in place of $(\Omega,I)$.

Thus, $\Omega'$ is in $I_a$. Also, (it appears that) b' is below the next element of $I_{a+1}$. The built-from-below requirement should ensure that comparison works correctly provided that we preserved enough structure for the passthrough ability. Here, that structure is the level of indiscernibility, and for the final level reached, the level of agreement with $\Omega$. Beyond that, built-from-below should let us get away with simulated structures: I' for b and (a different) $(\Omega',I')$ for c.

Starting with the minimum inner model L[U] with a measurable cardinal, the required indiscernibles can be constructed (using U and indiscernible in the model obtained by iterating U away) in a standard way. Also, it appears that we only care about indiscernibility of 1-tuples (with ordinal parameters). Furthermore, there is likely a weakening by using indiscernibles from a's iteration of the sharp that gets the exact strength of the notation system (perhaps for an extension of the system with additional terms). In that weakening, simulated indiscernibles (such as I') can be constructed in L (and other weak models) using a generic extension that makes Ω countable; L[($I_{a'}$:a'≤a)] is used in place of V. To get optimality, Ω (or the starting Ω) can be chosen as the least element of $I_{a+1}$ above d; the extent of I' will depend on b and (for Ω') c. Assuming it works and is optimal, the strength is that of "rudimentary set theory plus the sharp operation can be iterated any ordinal number of times".

REFERENCES

Arai, Toshiyasu (2015), "A simplified ordinal analysis of first-order reflection", arXiv:1506.05280.
Duchhardt, Christoph (2008), "Thinning Operators and $\Pi_4$-Reflection", thesis, http://d-nb.info/98.3428.3/34.
Jech, Thomas (2006). *Set Theory* (3rd edition). Springer.
Kanamori, Akihiro (2008). *The Higher Infinite: Large Cardinals in Set Theory from Their Beginnings* (2nd edition). Springer.
Rathjen, Michael (1994), "Proof theory of reflection", Annals of Pure and Applied Logic 68, pp. 181–224.
Rathjen, Michael (2006), "The art of ordinal analysis", International Congress of Mathematicians, II, Zürich, Eur. Math. Soc., pp. 45–69.
Rathjen, Michael (2014), "Omega-models and well-ordering principles", Foundational Adventures: Essays in Honor of Harvey M. Friedman (Neil Tennant (ed.)). (College Publications, London, 2014) 179-212.
Stegert, Jan-Carl (2010), "Ordinal proof theory of Kripke-Platek set theory augmented by strong reflection principles", thesis, https://nbn-resolving.org/urn:nbn:de:hbz:6-44449504436.
Taranovsky, Dmytro (2012), "Reflective Cardinals", arXiv:1203.2270.